\documentclass[12pt]{article}

\usepackage{amssymb}
\usepackage{amsmath}
\usepackage{color}
\usepackage{graphics}
\usepackage{epsfig}
\usepackage{hyperref}

\newtheorem{claim}{\bf \t}[part]

\newtheorem{Definition}{Definition}[part]

\newtheorem{Lemma}{Lemma}[part]

\newtheorem{Remark}{Remark}[part]
\newtheorem{Theorem}{Theorem}[part]

\numberwithin{Assumption}{section} \numberwithin{Corollary}{section}
\numberwithin{Definition}{section} \numberwithin{equation}{section}
\numberwithin{Example}{section} \numberwithin{Lemma}{section}
\numberwithin{Proposition}{section} \numberwithin{Remark}{section}
\numberwithin{Theorem}{section}

\def\v{\varepsilon}

\def\x{\xi}
\def\t{\theta}

\def\mb{\mathbf}
\def\a{\alpha}
\def\b{\beta}
\def\g{\gamma}
\def\d{\delta}
\def\l{\lambda}
\def\f{\frac}

\def\r{\rho}
\def\ra{\rightarrow}
\def\s{\sigma}
\def\z{\zeta}

\def\di{\displaystyle}
\def\i{\infty}

\def\text#1{{\rm #1}}
\begin{document}
\date{}
\title{\Large \bf Vacuum behaviors around rarefaction waves to 1D compressible Navier-Stokes equations with
density-dependent viscosity}
\author{\small \textbf{Quansen Jiu},$^{1,3}$\thanks{The research is partially
supported by National Natural Sciences Foundation of China (No.
10871133) and Project of Beijing Education Committee. E-mail:
qsjiumath@gmail.com}\quad
  \textbf{Yi Wang}$^{2,3}$\thanks{The research is partially
supported by National Natural Sciences Foundation of China (No.
10801128). E-mail: wangyi@amss.ac.cn.}\quad and \textbf{Zhouping
Xin}$^{3}$\thanks{The research is partially supported by Zheng Ge Ru
Funds, Hong Kong RGC Earmarked Research Grant CUHK4042/08P and
CUHK4040/06P, and a Focus Area Grant at The Chinese University of
Hong Kong. Email: zpxin@ims.cuhk.edu.hk}} \maketitle \small $^1$
School of Mathematical Sciences, Capital Normal University, Beijing
100048, China

\small $^2$ Institute of Applied Mathematics, AMSS, and Hua Loo-Keng
Key Laboratory of Mathematics, CAS, Beijing 100190, China

\small $^3$The Institute of Mathematical Sciences, Chinese
University of HongKong, HongKong\\

 {\bf Abstract:} In this paper, we study the large time asymptotic behavior
 toward rarefaction waves for solutions to
  the 1-dimensional compressible
Navier-Stokes equations with density-dependent viscosities for
general initial data whose far fields are connected by a rarefaction
wave to the corresponding Euler equations with one end state being
vacuum.
 First, a global-in-time weak solution around  the
rarefaction wave is constructed by approximating the system and
regularizing the initial data with general perturbations, and some a
priori uniform-in-time estimates for the energy and entropy are
obtained. Then it is shown that the density of any weak solution
satisfying the natural energy and entropy estimates will converge to
the rarefaction wave connected to vacuum with arbitrary strength in
super-norm time-asymptotically. Our results imply, in particular,
that the initial vacuum at far fields will remain for all the time
which are in contrast to the case of non-vacuum rarefaction waves
studied in \cite{JWX} where all the possible vacuum states will
vanish in finite time. Finally, it is proved that the weak solution
becomes regular away from the vacuum region of the rarefaction wave.



\section{Introduction } \setcounter{equation}{0}
\setcounter{Assumption}{0} \setcounter{Theorem}{0}
\setcounter{Proposition}{0} \setcounter{Corollary}{0}
\setcounter{Lemma}{0} In this paper, we consider the following
compressible and isentropic Navier-Stokes equations with
density-dependent viscosities
\begin{eqnarray}\label{(1.1)}
\begin{cases}
    \rho   _{t}+ (\rho u) _{x} = 0,   \qquad\qquad\qquad\qquad\qquad x\in\mathbf{R},~t>0,\cr
    (\rho u)_{t}+  \big {(} \rho u^2 +p(\r)\big {)}_x
    =(\mu(\r)  u_{x})_x,\cr
   \end{cases}
\end{eqnarray}
where $\rho (t, x)\geq 0 $, $u(t, x) $ represent the density and the
velocity of the gas, respectively. Let the pressure and viscosity
function be given by
\begin{eqnarray}
p(\r)=A \r^{\gamma },\qquad \mu(\r)=B\r^\a,~~\label{(1.2)}
\end{eqnarray}
respectively, where $\gamma >1$ denotes the adiabatic exponent,
$\a>0$ and $A,B>0$ are the gas constants. For simplicity, it is
assumed that $A=B=1$.

Consider the Cauchy problem (\ref{(1.1)}) with the initial values
\begin{equation}
  (\r,\r u)( 0, x)=(\r_0,m_0)(x)\rightarrow (\r_{\pm},m_{\pm}), ~~ \text{as} ~~ x \rightarrow
 \pm\infty,\label{(1.3)}
\end{equation}
where $\r_{\pm}\ge 0$,  $m_{\pm} $  are prescribed constants.

The large time  behavior of solutions to (\ref{(1.1)})-\eqref{(1.3)}
is expected to be closely related to that of the Riemann problem of
the corresponding Euler system
\begin{equation}
\left\{
\begin{array}{l}
\di \r_t+(\r u)_x=0,\\
\di (\r u)_t+(\r u^2+p(\r))_x=0,\\
\end{array}
\right. \label{(1.4)}
\end{equation}
with Riemann initial data
\begin{equation}
  (\r,\r u)( 0, x):=(\r_0^r,m^r_0)(x)=\left\{
  \begin{array}{ll}
(\r_-,m_-),~~~~~&x<0,\\
(\r_+,m_+),~~~~~&x>0.
  \end{array}
  \right.\label{Riemann}
\end{equation}

Different initial states (\ref{Riemann})  produce different type of
waves, namely, shock waves and rarefaction waves  for the
one-dimensional compressible isentropic Euler equations
\eqref{(1.4)}. However,  as pointed out by Liu-Smoller \cite{lius},
among the two nonlinear waves, i.e., shocks and rarefaction waves,
only rarefaction waves can be connected to vacuum. When vacuum
appears, the stability of rarefaction waves to the 1D compressible
Navier-Stokes equations is an important issue.

When the viscosity $\mu(\r)$ is a constant, there have been
extensive studies  on the stability of the rarefaction waves to the
1D compressible Navier-Stokes equations under the assumptions that
the rarefaction waves and the solutions  are away from the vacuum
(see \cite{FNP00}, \cite{JXZ}, \cite{KS}, \cite{L98},
\cite{Liu-Xin-1}, \cite{Matsumura-Nishihara-1},
\cite{Matsumura-Nishihara-2} and the references therein). However,
when vacuum appears, the well-known results by Hoff-Serre \cite{HD},
 Xin \cite{Xin} and Rozanova \cite{R} show that the solutions of the compressible
Navier-Stokes equations with constant viscosity may behave
singularly, in particular, in the case that the fluids jump to far
field vacuum. Liu, Xin and Yang first proposed in \cite{LXY} some
models of the compressible Navier-Stokes equations with
density-dependent viscosities to investigate the dynamics of the
vacuum. On the other hand, when deriving by Chapman-Enskog
expansions from the Boltzmann equation, the viscosity of the
compressible Navier-Stokes equations  depends on the temperature and
thus on the density for isentropic flows. Also, the viscous
Saint-Venant system for the shallow water, derived from the
incompressible Navier-Stokes equations with a moving free surface,
is expressed exactly as in \eqref{(1.1)}-\eqref{(1.2)} with $\a=1$
and $\g=2$ (see \cite{GP}). However,  there
 appear new mathematical challenges in dealing with such systems. In particular, these systems
become highly degenerate. The velocity  cannot even be defined in
the presence of vaccum and hence it is difficult to get  uniform
estimates for the velocity near vacuum. The global existence of
generak weak solutions to the compressible Navier-Stokes equations
with density-dependent viscosities or the viscous Saint-Venant
system for the shallow water model in the multi-dimensional case
remains open, and one can refere to \cite{BD3}, \cite{BDG},
\cite{GJX}, \cite{MV} for recent developments along this line.

There are  a large number of literatures on mathematical studies  of
\eqref{(1.1)}-\eqref{(1.2)} with various initial and boundary
conditions. If the initial density is assumed to be connected to
vacuum with discontinuities, Liu, Xin and Yang first obtained in
\cite{LXY} the local well-posedness of weak solutions. The global
well-posedness was obtained later by \cite{J}, \cite{JXZ},
\cite{OMM}, \cite{YYZ} respectively. The case of initial densities
connecting to vacuum continuously was studied by \cite{FZ},
\cite{VYZ}, \cite{YYZ} and \cite{YZ2} respectively. However, most of
these results concern with free boundary problems. Recently,
initial-boundary-value problems for the one-dimensional equations
\eqref{(1.1)}-\eqref{(1.2)} with $\mu(\rho)=\rho^\alpha
(\alpha>1/2)$ was studied by Li, Li and Xin
 in \cite{LLX} and the phenomena of vacuum vanishing and
blow-up of solutions were found there. The global existence of weak
solutions for the initial-boundary-value problems for spherically
symmetric compressible Navier-Stokes equations with
density-dependent viscosity was proved by Guo, Jiu and Xin in
\cite{GJX}. More recently, there are some results on global
existence of weak solutions to the Cauchy problem
\eqref{(1.1)}-\eqref{(1.3)}. The existence and uniqueness of global
strong solutions to the compressible Navier-Stokes equations
\eqref{(1.1)}-\eqref{(1.3)} were obtained by Mellet and Vasseur
\cite{MV} where no vacuum is permitted in the initial density and
for $0\leq\alpha<\f12$. However, the a priori estimates obtained in
\cite{MV} depend on the time interval thus do not give the
time-asymptotic behavior of the solutions. The first result about
the time-asymptotic behavior of the solutions to the Cauchy problem
\eqref{(1.1)}-\eqref{(1.3)} is obtained by Jiu-Xin \cite{JX}, where
the global existence and large time-asymptotic behavior of the weak
solutions were considered in the case that $\r_+=\r_-\geq0$ and
$u_+=u_-=0$. In the case that $\r_+=\r_->0$ and $u_+=u_-=0$, the
vanishing of the vacuum and the blow-up phenomena of the weak
solutions were also obtained in \cite{JX}. One of the key elements
in the analysis  in \cite{JX} is an interesting entropy estimate
which was observed first in \cite{Ka} for the one-dimensional case
and later established in \cite{BD1,BDL,BD2} for more general and
multi-dimensional cases due to the structure that the viscosity
coefficients vanish at the vacuum. This entropy estimate  provides
higher regularity of the density and played a crucial role  in
\cite{GJX,JX,LLX} for global existence and large time asymptotic
behaviors of  weak solutions.

The stability of rarefaction waves of the 1D compressible
Navier-Stokes equations with density-dependent viscosity was studied
in \cite{JWX} under general initial perturbations such that the
initial data and the solutions  may contain the vacuum. However, in
\cite{JWX}, the rarefaction wave itself is away from the vacuum. In
this paper, we are concerned with the
 the case when the rarefaction wave is
permitted to be connected to vacuum.

For definiteness,  we consider the case of  a 2-rarefaction wave
such that $\r_-=0, \r_+>0$ in (\ref{(1.3)}). Similar to \cite{JWX},
we will first construct a class of approximate solutions satisfying
some uniform estimates and furthermore prove the global existence of
weak solutions for the Cauchy problem of
\eqref{(1.1)}-\eqref{(1.3)}. To get the  uniform energy and entropy
estimates in time  to the approximate solutions, we combine the
elementary energy estimates and the  entropy estimates
 in an elaborate way. Note that the elementary energy
estimates and the entropy estimates  are coupled to each other due
to the underlying rarefaction wave. This is quite different from the
previous works on the global existence and the time-asymptotic
behavior of the solutions to Navier-Stokes equations \eqref{(1.1)}
with density-dependent viscosity where the elementary energy
estimates and the entropy estimates can be derived independently.
Moreover, compared with the case of non-vacuum rarefaction waves in
\cite{JWX}, some new  difficulties occur due to the degeneracies at
the vacuum states in the 2-rarefaction wave. To overcome these
difficulties, we first cut off the 2-rarefaction wave with vacuum
along the rarefaction wave curve and then derive some uniform
estimates with respect to both the approximations and the cut-off
process. More precisely, for any $\nu>0$, a suitably small
parameter, the cut-off 2-rarefaction wave will connect the state
$(\r,u)=(\nu,u_\nu)$ and $(\r_+,u_+)$ where $u_\nu$ can be obtained
explicitly and uniquely by the definition of the 2-rarefaction wave
curve. For any fixed $\nu>0$, one can obtain a weak solution to the
compressible Navier-Stokes equations \eqref{(1.1)}-\eqref{(1.3)}
with $(\r_-,m_-)$ replaced by $(\nu,\nu u_\nu)$ along the same line
as in our previous paper \cite{JWX}. Thus, in order to get the
solution to the original problem \eqref{(1.1)}-\eqref{(1.3)}, we
will derive some uniform estimates with respect to both the
approximations and the cut-off process. To this end, the
approximation parameters $\v$ and the cut-off parameter $\nu$ should
be chosen in an appropriate way. Thus, as a limit of this
approximate solution, a global weak solution to
\eqref{(1.1)}-\eqref{(1.3)} is shows to exist with the
uniform-in-time estimates \eqref{P12-1} and \eqref{P12-3}.

Next, we study the large-time asymptotic behavior of any weak
solutions to \eqref{(1.1)}-\eqref{(1.3)} under the uniform-in-time
bounds \eqref{P12-1} and \eqref{P12-3}. It is shown that
time-asymptotically, the density function tends to the rarefaction
wave connected to the vacuum in $L^\i$ norm. This time-asymptotic
behavior of the density function implies that the vacuum in the far
field is essential and  will maintain for all the time.  This is
quite different from the
 previous results in \cite{JWX} and \cite{LLX}  where all the possible vacuum states will
 vanish in finite time.
 At last, we prove that such a weak solution becomes regular away from
the vacuum
 region of the rarefaction wave by using the Di Giorgi-Moser
iteration and higher order energy estimates.

 \vskip 2mm
\noindent\emph{Notations.} Throughout this paper, positive generic
constants are denoted by $c$ and $C$, which are independent of $\v$,
$\nu$ and $T$, without confusion, and $C(\cdot)$ stands for some
generic constant(s) depending only on the quantity listed in the
parenthesis. For function spaces, $L^{p}(\Omega), 1\leq p\leq
\infty$, denote the usual Lebesgue spaces on $\Omega \subset
\mathbf{R}:=(-\infty,\infty)$. $W^{k,p}(\Omega)$ denotes the
$k^{th}$ order Sobolev space, $H^{k}(\Omega):=W^{k,2}(\Omega)$.

\section{Preliminaries and Main Results}

In this section we first describe the rarefaction wave connected to
the vacuum to the compressible Euler system \eqref{(1.4)}. Then an
approximate rarefaction wave will be constructed through the
Burger's equation and the main results of the paper will be given at
last.

\subsection{Rarefaction waves}

The Euler system (\ref{(1.4)}) is a strictly hyperbolic one for
 $\r>0$ whose characteristic fields are both genuinely nonlinear, that
 is, in the equivalent system
$$
\left(
\begin{array}{l}
\di \r\\
\di u
\end{array}
\right)_t + \left(
\begin{array}{cc}
\di u&\quad \r\\
\di p^\prime(\r)/\r&\quad u
\end{array}\right)\left(
\begin{array}{l}
\di \r\\
\di u
\end{array}\right)_x=0,
$$
the Jacobi matrix
$$
\left(
\begin{array}{cc}
\di u&\quad \r\\
\di p^\prime(\r)/\r&\quad u
\end{array}\right)
$$
has two distinct eigenvalues
$$
\l_1(\r,u)=u-\sqrt{p^\prime(\r)},\qquad
\l_2(\r,u)=u+\sqrt{p^\prime(\r)}
$$
with corresponding right eigenvectors
$$
r_i(\r,u)=(1,(-1)^i\f{\sqrt{p^\prime(\r)}}{\r})^t,\qquad i=1,2,
$$
such that
$$
r_i(\r,u)\cdot \nabla_{(\r,u)}\l_i(\r,u)=(-1)^i\f{\r
p^{\prime\prime}(\r)+2p^\prime(\r)}{2\r\sqrt{p^\prime(\r)}}\neq
0,\quad i=1,2, \quad {\rm if}~~\r>0.
$$
Define the $i-$Riemann invariant $(i=1,2)$ by
$$
\Sigma_i(\r,u)=u+(-1)^{i+1}\int^\r\f{\sqrt{p^\prime(s)}}{s}ds,
$$
such that
$$
\nabla_{(\r,u)}\Sigma_i(\r,u)\cdot r_i(\r,u)\equiv0,\qquad \forall
\r>0,u.
$$

 There are two families of rarefaction waves
to the Euler system (\ref{(1.4)})-\eqref{Riemann}. Here we consider
the case of $2-$rarefaction wave connected with vacuum, that is
$\r_-=m_-=0,\r_+>0$. Thus we can define the velocity at the positive
far field $u_+=\f{m_+}{\r_+}$. First we give the description of the
2-rarefaction wave connected with vacuum, see also in details in
\cite{lius}. From the fact that $2-$Riemann invariant is constant:
$$
\Sigma_2(\r_-=0,u_-)=\Sigma_2(\r_+,u_+),
$$
we can define the velocity $u_-$ which is the speed of the gas
coming into the vacuum region. Then the entropy condition
$\l_2(\r_-=0,u_-)<\l_2(\r_+,u_+)$ is always satisfied. This
$2-$rarefaction wave connecting the vacuum $\r_-=0$ to $(\r_+,u_+)$
is the self-similar solution $(\r^r,u^r)(\x),~(\x=\f xt)$ of
(\ref{(1.4)})-\eqref{Riemann} defined by
\begin{equation}
\begin{array}{l}
\qquad\qquad\qquad\qquad~~ \r^r(\x)=0, ~\quad {\rm if}~~\x<\l_2(0,u_-)=u_-,\\
\l_2(\r^r(\x),u^r(\x))=\left\{
\begin{array}{ll}
\di \x, &\di {\rm if}~~u_-\leq\x\leq\l_2(\r_+,u_+),\\
\di \l_2(\r_+,u_+), &\di {\rm if}~~\x>\l_2(\r_+,u_+),
\end{array} \right.
\end{array}\label{2-R1}
\end{equation}
and
\begin{equation}
\Sigma_2(\r^r(\x),u^r(\x))=\Sigma_2(\r_+,u_+)=\Sigma_2(0,u_-).\label{2-R2}
\end{equation}
Thus we can define the momentum of 2-rarefaction wave by
\begin{equation}m^{r}(\x)=\left\{\begin{array}{ll}\r^{r}(\x) u^{r}(\x),\qquad &{\rm if}~~~\r^r(\x)>0,\\
0,&{\rm if}~~~ \r^r(\x)=0.
\end{array}\right.\label{2-R3}
\end{equation}
In this paper, we consider the time-asymptotic behavior toward such
 rarefaction waves of solutions to the compressible Navier-Stokes
equations \eqref{(1.1)} with density-dependent viscosities.

\subsection{Approximate rarefaction waves}
Consider the Riemann problem for the inviscid Burgers equation:
\begin{align}\label{bur}
\left\{\begin{array}{ll}
w_t+ww_x=0,\quad \,t>0, ~ x\in\mathbf{R}, \\
w(x,0)=\left\{\begin{array}{ll}
w_-,&x<0,\\
w_+,&x>0.
\end{array}
\right.
\end{array}
\right.
\end{align}
If $w_-<w_+$, then the Riemann problem $(\ref {bur})$ admits a
rarefaction wave solution $w^r(x, t) = w^r(\f xt)$ given by
\begin{align}\label{abur}
w^r(\f xt)=\left\{\begin{array}{lr}
w_-,&\f xt\leq w_-,\\
\f xt,&w_-\leq \f xt\leq w_+,\\
 w_+,&\f xt\geq w_+.
\end{array}
\right.
\end{align}

Consider   the solution to the following Cauchy problem for Burgers
equation
\begin{eqnarray}\label{(2.1)}
\begin{cases}
 \di w_{t}+ww_{x}=0,  \quad \,t>0, ~ x\in\mathbf{R}, \cr
 \di w( 0,x):=w_0( x) =\f{w_++w_-}{2}+\f{w_+-w_-}{2}K_q
\int^{\eta x }_0(1+y^2)^{-q}\,dy.
\end{cases}
\end{eqnarray}
Here $q\geq 2$ is some fixed constant, and $K_q$ is a constant such
that $\di K_q\int^{\infty}_0(1+y^2)^{-q} dy=1 $, and $\eta$ is a
small positive constant to be determined later. It is easy to see
that the solution to this problem is given by
\begin{eqnarray}
w(t,x)=w_0(x_0(t, x)),\quad\quad x=x_0(t, x)+w_0(x_0(t, x))t.
\end{eqnarray}
Then the following properties hold (see
\cite{Matsumura-Nishihara-2}).

\begin{Lemma}  Let $w_-<w_+$, the Cauchy problem
$\eqref{(2.1)}$ has a unique smooth solution $w(t,x)$ satisfying

i) ~ $w_-< w(t,x)<w_+,~w_x(t,x)> 0 $;

ii) ~ For any $p$ $(1\leq p\leq \infty)$, there exists a constant
$C_{pq}$ such that
\begin{eqnarray}\nonumber
&&\|w(t,\cdot)-w^r(\f \cdot t)\|_{L^p(\mathbf{R})}\leq
 C_p \d_r\eta^{-\f1p},\cr
&& \parallel w_x(t)\parallel_{L^p(\mathbf{R})}\leq
C_{pq}\min\big{\{}\delta_r\eta^{1- \f1p},~
\delta_r^{\f1p}t^{-1+\f1p}\big{\}},   \cr &&\parallel
w_{xx}(t)\parallel_{L^p(\mathbf{R})}\leq
C_{pq}\min\big{\{}\delta_r\eta^{2-\f1p},~
\eta^{(1-\f{1}{2q})(1-\f1p)}\d_r^{-\f{(p-1)}{2pq}}t^{-1-\f{(p-1)}{2pq}}\big{\}},
\end{eqnarray}
where $\d_r=w_+-w_-,$ and $C_p,C_{pq}$ are independent of
 $t$;

 iii) ~ $\sup\limits_{x\in\mathbf{R}} |w(t,x)-w^r(\f
xt)|\rightarrow 0$, as $t\rightarrow \infty$.
\end{Lemma}

We now turn to rarefaction waves to the Euler system
\eqref{(1.4)}-\eqref{Riemann}. Set $\l_2(\r_\pm,u_\pm)=w_\pm$ with
$\r_-=0$ in \eqref{bur}. Then the unique solution $(\r^r,u^r)(\x)$
in \eqref{2-R1}-\eqref{2-R3} to the Riemann problem
\eqref{(1.4)}-\eqref{Riemann} can also be expressed in terms of
$w^r(\f xt)$ in \eqref{abur}, by
\begin{equation}
\begin{array}{l}
\di \l_2(\r^r(\f xt),u^r(\f xt))=w^r(\f xt),\\
\di \Sigma_2(\r^r(\f xt),u^r(\f xt))=\Sigma_2(\r_{\pm},u_\pm).
\end{array}\label{R2}
\end{equation}
Correspondingly,  an approximate $2-$rarefaction wave $(\bar\r, \bar
u) (t,x)$ can be defined by
\begin{equation}
\begin{array}{l}
\di \l_2(\bar\r(t,x),\bar u(t,x))=w(1+t,x),\\
\di \Sigma_2(\bar\r(t,x),\bar u(t,x))=\Sigma_2(\r_{\pm},u_\pm).
\end{array}\label{a2R}
\end{equation}
It can be checked that $(\bar\r, \bar u) (t,x)$ also satisfies the
Euler system
\begin{equation}\label{(2.4)}
\left\{
\begin{array}{l}
\di \bar\r_t+(\bar\r\bar u)_x=0\\
\di (\bar\r \bar u)_t+(\bar\r\bar u^2+p(\bar\r))_x=0,
\end{array}
\right.
\end{equation}
and properties listed in the following Lemma.

\begin{Lemma}\label{rarefaction-1} The approximate $2-$rarefaction wave $(\bar\r, \bar u) (t,x)$
defined in \eqref{a2R} satisfies
\begin{itemize}
\item[(1)]  $\bar\r_x>0,\quad \bar u_x>0,\quad \bar u_x=\sqrt \g\bar\r^{\f{\g-3}2}\bar\r_x$;
\item[(2)]  For any $p$ ($1\leq p\leq \infty$), there exists a constant
$C_{pq}$ such that
$$
\begin{array}{ll}
\di \|(\bar\r,\bar u)(t,\cdot)-(\r^r,u^r)(\f \cdot
t)\|_{L^p(\mathbf{R})}\leq
 C_p(w_+-w_-)\eta^{-\f1p},\\
\di  \|\bar u_x(t,\cdot)\|_{L^p(\mathbf{R})}\leq
C_{pq}\min\{\d\eta^{1-\f1p},
\d^\f1p(1+t)^{-1+\f1p}\} ,\\
\di\|\bar u_{xx}(t,\cdot)\|_{L^p(\mathbf{R})}\leq
C_{pq}\min\{\d\eta^{2-\f1p},
\eta^{(1-\f{1}{2q})(1-\f1p)}\d^{-\f{p-1}{2pq}}(1+t)^{-1-\f{p-1}{2pq}}+\d^{\f1p}(1+t)^{-2+\f1p}\},
 \end{array}$$
 where $\d=|\r_+-\r_-|+|u_+-u_-|$ is the strength of the rarefaction
 wave;
 \item[(3)]   $ \lim\limits_{t\rightarrow
\infty}\sup\limits_{x\in\mathbf{R}}\big| \bar\r (t,x)-\r^{r}(\f
x{t})\big| =0$.
\end{itemize}
\end{Lemma}

{\begin{Remark} For any $1<p\leq +\i$,
$$
\int_0^T\|\bar u_{xx}(t,\cdot)\|_{L^p(\mathbf{R})}dt\leq C,
$$
where $C$ is independent of $T$. Note that in the case $p=1$, the
constant $C$ in the above estimates is not uniform in $T$. Moreover,
the following estimate holds:
$$
\int_0^T\|\bar u_{xx}(t,\cdot)\|_{L^\i(\mathbf{R})}dt\leq
C\eta^{\f2{4q+1}}\int_0^T(1+t)^{-1-\f1{4q+1}}dt\leq
C\eta^{\f2{4q+1}},
$$
where $C$ is independent of $T$.
\end{Remark}

\subsection{Main Results}
Set
\begin{equation}\label{10-21-3}
\begin{array}{ll}
\di \Psi(\r,\bar\r)&\di =\int_{\bar\r}^\r\f{p(s)-p(\bar\r)}{s^2}ds\\
&\di
=\f{1}{(\g-1)\r}\Big[\r^\g-\bar\r^\g-\g\bar\r^{\g-1}(\r-\bar\r)\Big].
\end{array}
\end{equation}
The initial values are assumed to satisfy:
 \begin{equation}
 \left\{
 \begin{array}{ll}
  &\rho_0\ge 0; \ \ \ m_0=0\ \ a.e.~ {\rm on}\ \{x\in \mathbf{R}|\rho_0(x)=0\};\\
 &
 \quad (\rho_0^{\a-\f12})_x\in
 L^2(\mathbf{R}),\quad \r_0\Psi(\r_0,\r_\pm)\in L^1(\mathbf{R}^\pm)~~ {\rm with}~~ \r_-=0;\\
 &\di \r_0(\frac{m_0}{\rho_0}-\frac{m_\pm}{\rho_\pm})^2\in L^1(\mathbf{R}^{\pm}),\quad \r_0(\frac{m_0}{\r_0}-\frac{m_\pm}{\rho_\pm})^{3}\in
 L^1(\mathbf{R}^\pm),
 \end{array}
\right.\label{id}
 \end{equation}
where  $u_-:=\f{m_-}{\r_-}$ since $\r_-=m_-=0$. Note that \eqref{id}
implies that $\rho_0\in C_B(\mathbf{R})$ which is the space of
bounded and continuous functions.

Equivalently, the assumptions \eqref{id} can be rewritten as
\begin{equation}
 \left\{
 \begin{array}{ll}
  &\rho_0\ge 0; \ \ \ m_0=0\ \ a.e.~ {\rm on}\ \{x\in \mathbf{R}|\rho_0(x)=0\};\\
 &
 \quad (\rho_0^{\alpha-\frac 12})_x\in
 L^2(\mathbf{R}),\quad \r_0\Psi(\r_0,\bar\r_0)\in L^1(\mathbf{R});\\
 &\di \r_0(\frac{m_0}{\rho_0}-\bar u_0)^2\in L^1(\mathbf{R}),\quad \r_0(\frac{m_0}{\r_0}-\bar u_0)^{3}\in
 L^1(\mathbf{R}),
 \end{array}
\right.\label{(2.5)}
 \end{equation}
 where $(\bar\r_0,\bar u_0):=(\bar\r,\bar u)(0,x)$ is
the initial values of the approximate $2-$rarefaction wave
$(\bar\r,\bar u)(t,x)$ constructed in \eqref{a2R}.

Before stating the main results, we give the definition of  weak
solutions to \eqref{(1.1)}-\eqref{(1.3)} with the far fields
$(\r_\pm,m_\pm)$. Let $T>0$ be given. For any far fields
$(\r_\pm,u_\pm)$ satisfying $\r_\pm\ge 0$ and any smooth  functions
$(\tilde\r,\tilde u)(t,x)$ connecting with them, we define

\begin{Definition}\label{def} A pair $(\rho, u)$ is said to be a weak solution
to \eqref{(1.1)}-\eqref{(1.3)} with the far fields $(\r_\pm,u_\pm)$
provided that there exists a smooth functions $(\tilde\r,\tilde
u)(t,x)$  with the same far fields $(\r_\pm,u_\pm)$ and
$\tilde\r\geq0$, such that

(1) $\rho\geq 0$ a.e., and
\begin{eqnarray*}
&&\rho\in L^\infty(0,T; L^\infty(\mathbf{R})))\cap C([0,T]\times\mathbf{R}),\quad \r\Psi(\r,\tilde\r)\in L^\i(0,T;L^1(\mathbf{R})),\\
&& (\rho^{\a-\f12})_x\in
L^\infty(0,T;L^2(\mathbf{R})),\qquad\qquad~~ \sqrt{\rho}(u-\tilde
u)\in L^\infty(0,T;L^2(\mathbf{R}));
\end{eqnarray*}

 (2) For any $t_2\geq
t_1\ge 0$ and any $\zeta\in C_0^1({\mathbf{R}}\times[t_1,t_2])$, the
mass equation \eqref{(1.1)} holds in the following sense:
\begin{eqnarray}
\int_{\mathbf{R}} \rho\zeta
dx|_{t_1}^{t_2}=\int_{t_1}^{t_2}\int_{\mathbf{R}}(\rho
\zeta_t+\sqrt\rho \sqrt\rho u\zeta_x)dxdt; \label{mass}
\end{eqnarray}

(3) For any $\psi\in C_0^\infty({\mathbf{R}}\times(0,T))$, it holds
that
\begin{eqnarray}
&&\int_0^T\int_{\mathbf{R}}\Big\{\sqrt{\rho}\sqrt{\rho}u\psi_t+\big[(\sqrt{\rho}u)^2+\rho^\gamma\big]\psi_x\Big\}dxdt+\langle\rho^\a
u_x, \psi_x\rangle=0, \label{momentum}
\end{eqnarray}
where the diffusion term makes sense if written as
\begin{equation}
\begin{array}{ll}
\di \langle\rho^\a u_x, \psi_x\rangle&\di
=-\int_0^T\int_{\mathbf{R}} \rho^{\a-\f12}\sqrt{\rho} u\psi_{xx}
dxdt-\f{2\a}{2\a-1}\int_0^T\int_{\mathbf{R}}
(\rho^{\a-\f12})_x\sqrt{\rho} u\psi_x dxdt. \end{array}
\label{diffusion}
\end{equation}
\end{Definition}

The first main result in this paper reads as

\begin{Theorem}\label{w-e}Let $\alpha$ and $\gamma$ satisfy that
\begin{equation}
1<\g\leq2,\quad{\rm and}\quad 1\leq\a\leq\f{\g+1}{2},\label{alpha}
\end{equation} and suppose that
\eqref{(2.5)} holds. Then the Cauchy problem
\eqref{(1.1)}-\eqref{(1.3)} admits a global weak solution
$(\rho(x,t), u(x,t))$ in the sense of Definition \ref{def} with the
smooth function $(\tilde\r,\tilde u)$ replaced by the approximate
rarefaction wave $(\bar\r,\bar u)$. Furthermore, this weak solution
$(\rho(x,t), u(x,t))$ satisfies
\begin{equation}
\begin{array}{l}
\di\rho\ge 0, \quad \max_{(x,t)\in \mathbf{R}\times[0,T]}\rho\le
C,\qquad\rho\in  C(\mathbf{R}\times
[0,T]),\qquad\qquad\qquad\qquad\qquad\qquad\quad
\end{array}
\label{P12-1}
\end{equation}
\begin{equation}
\begin{array}{ll}
\di \sup_{t\in[0,T]}\int_\mathbf{R} \big[|\sqrt{\rho}(u-\bar
u)|^2+(\rho^{\alpha-\frac12})_x^2+\r\Psi(\r,\bar\r)\big]
dx\\
\di +\int_0^T\int_\mathbf{R}
\big\{[(\rho^{\frac{\gamma+\alpha-1}{2}}-\bar\rho^{\frac{\gamma+\alpha-1}{2}})_x]^2+\bar
u_x\r\Psi(\r,\bar\r)+\r\bar u_x(u-\bar u)^2+\Lambda(x,t)^2\big\}
dxdt\le C,
\end{array}\label{P12-3}
\end{equation}
where $C$ is an absolute constant depending on the initial data but
independent of $T$ and $\Lambda(x,t)\in L^2(\mb{R}\times (0,T))$
satisfies
\begin{eqnarray}
&&\int_0^T\int_\mathbf{R}\r^{\f\a2} \Lambda\varphi
dxdxt=-\int_0^T\int_\mathbf{R} \rho^{\alpha-\frac12}\sqrt{\rho}
(u-\bar u) \varphi_x
dxdt\nonumber\\
&&-\frac{2\alpha}{2\alpha-1}\int_0^T\int_\mathbf{R}
(\rho^{\alpha-\frac12})_x\sqrt{\rho} (u-\bar u)\varphi dxdt,\quad
\forall \varphi\in C_0^\i(\mb{R}\times (0,T)). \label{diffusion1}
\end{eqnarray}
\end{Theorem}

\begin{Remark} It should be noted that there is no
requirement on the sizes of the strength of the rarefaction wave and
the perturbations. The class of initial perturbations given by
\eqref {(2.5)} is quite large compared with those for the case of
constant viscosities, \cite{Liu-Xin-1},
\cite{Matsumura-Nishihara-1}, \cite{Matsumura-Nishihara-2}.
\end{Remark}

\begin{Remark} The important case of the shallow water
model, i.e., $\a=1,\g=2$, is included in our theorem.
\end{Remark}

\begin{Remark} For any far fields $(\r_\pm,u_\pm)$ satisfying $\r_\pm\ge 0$ and any smooth functions $(\tilde\r,\tilde
u)(t,x)$  connecting with them, one can also obtain the existence of
weak solutions in the sense of Definition \ref{def} in a similar way
(see \cite{Li-J} for the weak solutions in the case of $\r_\pm>0$
and $u_\pm=0$). However, in order to get the uniform in time
estimates in \eqref{P12-1}-\eqref{P12-3}, it seems that the far
fields $(\r_\pm,u_\pm)$ should be specified and in Theorem \ref{w-e}
the smooth function $(\tilde\r,\tilde u)(t,x)$ is replaced by the
approximate rarefaction wave $(\bar\r,\bar u).$
\end{Remark}

The next result concerns  on the asymptotic behaviors of the weak
solution, which can be stated as

\begin{Theorem}\label{asymptotic-b}
 Let $\a$ and $\g$ satisfy \eqref{alpha} and suppose that \eqref{id} holds. Suppose that  $(\rho, u)(x,t)$ is
  a global weak
solution of the Cauchy problem \eqref{(1.1)}-\eqref{(1.3)} in the
sense of Definition \ref{def} satisfying  \eqref{P12-1} and
\eqref{P12-3}. Then it holds that
\begin{equation}
\lim_{t\rightarrow+\i}\sup_{x\in\mathbf{R}}|\r(t,x)-\bar\r(t,x)|=0.\label{con}
\end{equation}
Consequently,
\begin{eqnarray}
\lim_{t\to\infty}\sup_{x\in \mathbf{R}}|\rho(x,t)-\rho^r(\f xt)|=0.
\label{th14-1}
\end{eqnarray}
\end{Theorem}

\begin{Remark}
A direct consequence of Lemma \ref{rarefaction-1} (2), \eqref{P12-3}
and \eqref{th14-1} is the following time-asymptotic behavior of the
density function:
\begin{eqnarray}
\lim_{t\to\infty}\|\rho(\cdot,t)-\rho^r(\f \cdot
t)\|_{L^p(\mathbf{R})}=0,\quad \forall~ 2<p\leq +\i. \label{th14-2}
\end{eqnarray}
\end{Remark}

\begin{Remark}
Theorem \ref{asymptotic-b} implies that for any weak solution
$(\rho,u)$ to the Cauchy problem \eqref{(1.1)}-\eqref{(1.3)} with
the far fields given by the vacuum state and $(\r_+,u_+)$, if
$(\rho,u)$ satisfies the bounds \eqref{P12-1} and \eqref{P12-3},
then the density function converges to the 2-rarefaction wave to the
corresponding Euler equations connecting the vacuum state and
$(\r_+,u_+)$ in sup-norm as $t$ tends to infinity. Consequently, the
initial vacuum at far field will remain for all the time, which is
contrast to the case of non-vacuum rarefaction waves studied in
\cite{JWX} where all the possible vacuum states will vanish.
\end{Remark}

Finally, we can obtain  the following higher regularity to the
velocity function $u(t,x)$ to a global weak solution $(\r,u)(t,x)$
of the Cauchy problem \eqref{(1.1)}-\eqref{(1.3)} in the sense of
Definition \ref{def} satisfying  \eqref{P12-1} and \eqref{P12-3} in
the region away from the vacuum region of 2-rarefaction wave
$(\r^r,u^r)(\x)$.

\begin{Theorem}\label{regularity} (Regularity of the solution away from the
vacuum) Let $(\rho, u)$ be a weak solution to  the Cauchy problem
\eqref{(1.1)}-\eqref{(1.3)} satisfying \eqref{P12-1} and
\eqref{P12-3}. For any fixed $\sigma>0$, there exist a straight line
$x=\l_2^\sigma t$ with
$\l_2^\sigma=\l_2(\r,u)|_{(\r,u)=(\sigma,u_\sigma)}$ defined in
\eqref{sigma} and a large time $T_{\sigma}$, such that if $
(t,x)\in\Omega_{\sigma}:=\{(t,x)|t>T_{\sigma}, x>\l^\sigma_{2}t\},$
then the density has the lower bound
$$
\r(t,x)\geq \f\s2.
$$
 Furthermore, for any
$(t_*,x_*)\in \Omega_{\sigma}$ and for any $r,s>0$ such that
$Q_{r,s}^*:=B_r(x_*)\times(t_*,t_*+s] \subset\Omega_{\sigma}$ with
$B_r(x_*)$ being the ball with the radius $r$ and the center $x_*$,
there exists a constant $\a_0\in(0,1),$ such that
\begin{equation}
\begin{array}{ll}
\di u\in C^{\alpha_0,\f{\alpha_0}{2}}_{\rm loc}(Q_{r,s}^*),\quad
u\in
L^\i_{\rm loc}(t_*,t_*+s,H^1_{\rm loc}(B_r(x_*))),\\
\di u_t\in L^2_{\rm loc}(Q_{r,s}^*),\quad u\in L^2_{\rm
loc}(t_*,t_*+s,H^2_{\rm loc}(B_r(x_*))).
\end{array}
\end{equation}
\end{Theorem}

\section{Existence of a weak solution}

We first study the following approximate system:
\begin{eqnarray}\label{(3.1)}
\begin{cases}
    \rho   _{t}+ (\rho u) _{x} = 0,   \qquad\qquad\qquad\qquad\qquad\qquad x\in\mathbf{R},~~ t>0,\cr
    (\rho u)_{t}+  \big {(} \rho u^2 +p(\r)\big {)}_x
    =(\mu_\v(\r)  u_{x})_x,\cr
   \end{cases}
\end{eqnarray}
where $\mu_\v(\r)=\r^\a+\v\r^\t$ with $\v>0$ and $\t=\f12$. This
kind of  the approximation $\v\r^\t$ was first used in \cite{JXZ}
and $\t=\f12$ is crucial in getting the lower bound of the
approximation density function to enure the existence of approximate
solutions.

To overcome the difficulty caused by the vacuum in the rarefaction
wave, we first cut off the rarefaction wave  along the wave curve.
More precisely, for any $\nu>0$ suitably small and to be determined,
let $(\nu,u(\nu))$ be the state  such that
$$
\Sigma_2(\nu,u(\nu))=\Sigma_2(\r_+,u_+),
$$
where 2-Riemann invariant $ \Sigma_2(\r,u)=u-\f{2\sqrt\g}{\g-1}
\r^{\f{\g-1}{2}}.$ Then $(\nu,u(\nu))$ is connected to $(\rho_+,
u_+)$ by a non-vacuum 2-rarefaction wave  given by
$(\r^{r}_{\nu},u^{r}_{\nu})(\f xt)$. Then it holds that
$$
|(\r^r,m^r)(\f xt)-(\r^{r}_{\nu},m^{r}_{\nu})(\f xt)|\leq C\nu.
$$
One can compute that $u(\nu)=\f{2\sqrt\g}{\g-1}
\nu^{\f{\g-1}{2}}+\Sigma_2(\r_+,u_+).$ So the corresponding smooth
approximate $2-$rarefaction wave $(\bar\r_\nu,\bar u_\nu)$ described
in Section 2.2 can be constructed by setting
$$
\l_2(\nu,u(\nu))=w_-,\qquad \l_2(\r_+,u_+)=w_+.
$$
Consequently, $(\bar\r_\nu,\bar u_\nu)$ will converge to
$(\bar\r,\bar u)$ point-wisely as $\nu$ tends to zero. In fact,
$\nu=\nu(\v)$ will be chosen suitably such that $\nu(\v)\rightarrow
0$ as $\v$ tends to zero.

The initial values $(\r_{0},m_{0})$ can be regularized in a similar
way as in \cite{JWX} such that
 \begin{equation} (\r,\r
u)(0,x)=(\r_{0\v,\nu},m_{0\v,\nu})\rightarrow\left\{\begin{array}{ll}
(\nu,\nu u(\nu)), ~~ \text{as} ~~ x \rightarrow
 -\infty,\\
 (\r_+,m_+), ~~
\text{as} ~~ x \rightarrow
 +\infty,\end{array}
 \right.\label{10-21-1}
\end{equation}
and
\begin{eqnarray}
\r_{0\v,\nu}(x)\geq \min\{\nu,\f{1}{2}\v^{\f2{2\a-1}}\},~~~\forall
x\in \mathbf{R},
\end{eqnarray}
for suitably small $\v,\nu>0$.

Furthermore,  $\r_{0\v,\nu}$ satisfies
$$
\r_{0\v,\nu}\Psi(\r_{0\v,\nu},\bar\r_{0\nu})\rightarrow
\r_0\Psi(\r_{0},\bar\r_0) ~{\rm in}~L^1(\mathbf{R}),\quad
(\rho_{0\v,\nu}^{\alpha-1/2})_x\to (\rho_{0}^{\alpha-1/2})_x \ {\rm
in}\ L^2(\mathbf{R}).
$$

Since in the following paper $\a, \g$ satisfy \eqref{alpha} and we
choose $\nu=\v^{\f23}$, it holds that
$$
\int_\mathbf{R} \v^2[(\ln\r_{0\v,\nu})_x]^2 dx
=(\f{\v}{\a-\f12})^2\int_\mathbf{R}
\r_{0\v,\nu}^{1-2\a}[(\r_{0\v,\nu}^{\a-\f12})_x]^2 dx\leq C.
$$

While $m_{0\v,\nu}$ satisfies
$$
\r_{0\v,\nu}(\f{m_{0\v,\nu}}{\r_{0\v,\nu}}-\bar
u_{0\nu})^{2}\rightarrow\r_0(\f{m_0}{\r_0}-\bar u_0)^{2}~~{\rm in}~~
L^1(\mathbf{R}),
$$
and
$$
\r_{0\v,\nu}(\f{m_{0\v,\nu}}{\r_{0\v,\nu}}-\bar
u_{0\nu})^{3}\rightarrow\r_0(\f{m_0}{\r_0}-\bar u_0)^{3}~~{\rm in}~~
L^1(\mathbf{R}).
$$

For any fixed $T>0$ and  $\v,\nu>0$, we will first construct smooth
approximate solutions $(\r_{\v,\nu},u_{\v,\nu})(x,t)$ to
(\ref{(3.1)})  with initial values $(\r,\r
u)(0,x)=(\r_{0\v,\nu},m_{0\v,\nu})(x)$ defined in $[0,T]$. To do
this, a key step is to get the lower bound of the density. Then the
global existence of  weak solutions to \eqref{(1.1)}-\eqref{(1.3)}
can be proved by compactness arguments. We intend to deduce the
uniform energy and entropy estimates with respect to $\v,\nu$ such
that one can pass to the limit $\v,\nu\rightarrow0$. Due to the
closeness to the vacuum of the rarefaction wave, we will have to
combine suitably the elementary energy estimates with the entropy
estimates to get the following estimates which are crucial to prove
our main results.

\begin{Lemma}\label{lemma1} Suppose that the conditions in Theorem
\ref{w-e} are satisfied and that $(\rho_{\v,\nu}, u_{\v,\nu})$ is a
smooth solution to \eqref{(3.1)} satisfying $\rho_{\v,\nu}>0$. Then
for any $T>0$ and $\v,\nu>0$ with $\v\ln(1+T)\leq C$ and $\v,\nu$
sufficiently small, the following estimate holds
\begin{equation}
\begin{array}{ll}
\di
\sup_{t\in[0,T]}\int_{\mathbf{R}}\Big\{\r_{\v,\nu}(u_{\v,\nu}-\bar
u_\nu)^2+\Big[\big(\r_{\v,\nu}^{\a-\f12}\big)_x\Big]^2+\v^2\big[(\ln\r_{\v,\nu})_x\big]^2+\r_{\v,\nu}\Psi(\r_{\v,\nu},\bar\r_\nu)\Big\}(x,t)dx
\\
\di+\int_0^T\int_{\mathbf{R}}\bigg\{(\bar u_{\nu})_x\Big[p(\r_{\v,\nu})-p(\bar\r_\nu)-p^\prime(\bar\r_\nu)(\r_{\v,\nu}-\bar\r_\nu)\Big]+\r_{\v,\nu} (\bar u_{\nu})_x(u_{\v,\nu}-\bar u_\nu)^2\\
\di \qquad +\r_{\v,\nu}^\a\Big[(u_{\v,\nu}-\bar
u_\nu)_x\Big]^2+\Big[(\r^{\frac{\gamma+\alpha-1}{2}}_{\v,\nu}-\bar\r^{\frac{\gamma+\alpha-1}{2}}_\nu)_x\Big]^2
\bigg\}(x,t)dxdt\leq C.
\end{array}\label{(3.3)}
\end{equation}
where $C>0$ is a universal  constant  independent of $\v,\nu$ and
$T$. \end{Lemma}

In the following, the subscripts $\v,\nu$ in the approximate
solution $(\r_{\v,\nu},u_{\v,\nu})(t,x)$ and  the subscripts $\nu$
in the approximate solution $(\bar\r_{\nu},\bar u_{\nu})(t,x)$ will
be omitted for simplicity.

 {\em
Proof:} {\it Step 1.  Energy Equality}

 It follows from $(\ref{(3.1)})_2$ that
\begin{equation}
\r u_t+\r uu_x+p(\r)_x=(\mu_\v(\r)u_x)_x.\label{10-21-2}
\end{equation}
Subtracting \eqref{10-21-2} from the second equation of
\eqref{(2.4)} gives
\begin{equation}
\r(u-\bar u)_t+\r u(u-\bar u)_x+(p(\r)-p(\bar\r))_x+(\r-\bar\r)\bar
u_t+(\r u-\bar\r\bar u)\bar u_x=(\mu_\v(\r)(u-\bar
u)_x)_x+(\mu_\v(\r)\bar u_x)_x.\label{(3.5)}
\end{equation}
Multiplying \eqref{(3.5)} by $u-\bar u$ yields
\begin{equation}
\begin{array}{ll}
\di \Big[\f{\r(u-\bar u)^2}{2}\Big]_t+\Big[\f{\r u(u-\bar
u)^2}{2}\Big]_x+(u-\bar u)(p(\r)-p(\bar\r))_x-\Big[\mu_\v(\r)(u-\bar
u)(u-\bar u)_x\Big]_x\\[2mm]
\di +\mu_\v(\r)\Big[(u-\bar u)_x\Big]^2=\Big[\mu_\v(\r)\Big]_x\bar
u_x(u-\bar u)+\mu_\v(\r)\bar u_{xx}(u-\bar u)-\Big[(\r-\bar\r)\bar
u_t+(\r u-\bar\r\bar u)\bar u_x\Big](u-\bar u).
\end{array}\label{(3.6)}
\end{equation}
Note that $\Psi(\r,\bar\r)$ defined in \eqref{10-21-3} satisfies
\begin{equation}
\begin{array}{l}
\di \Big[\r\Psi(\r,\bar\r)\Big]_t+\Big[\r
u\Psi(\r,\bar\r)\Big]_x+(u-\bar u)_x(p(\r)-p(\bar\r))+\bar
u_x\Big[\r^\g-\bar\r^\g-\g\bar\r^{\g-1}(\r-\bar\r)\Big]\\
\di =-\f{p(\bar\r)_x}{\bar\r}(\r-\bar\r)(u-\bar u).
\end{array}\label{(3.7)}
\end{equation}
It follows from \eqref{(3.6)} and \eqref{(3.7)} that
\begin{equation}
\begin{array}{ll}
\di \Big[\f{\r(u-\bar u)^2}{2}+\r\Psi(\r,\bar\r)\Big]_t+H_{1x}(t,x)
+\mu_\v(\r)\Big[(u-\bar u)_x\Big]^2+\bar
u_x\Big[\r^\g-\bar\r^\g-\g\bar\r^{\g-1}(\r-\bar\r)\Big]\\\di
=\Big[\mu_\v(\r)\Big]_x\bar u_x(u-\bar u)+\mu_\v(\r)\bar
u_{xx}(u-\bar u)-\Big[(\r-\bar\r)\bar u_t+(\r u-\bar\r\bar u)\bar
u_x+\f{p(\bar\r)_x}{\bar\r}(\r-\bar\r)\Big](u-\bar u),
\end{array}\label{(3.8)}
\end{equation}
where
$$
H_1(t,x)=\f{\r u(u-\bar u)^2}{2}+\r u\Psi(\r,\bar\r)+(u-\bar
u)(p(\r)-p(\bar\r))-\mu_\v(\r)(u-\bar u)(u-\bar u)_x.
$$
Since
$$
(\r-\bar\r)\bar u_t+(\r u-\bar\r\bar u)\bar
u_x+\f{p(\bar\r)_x}{\bar\r}(\r-\bar\r)=\r(u-\bar u)\bar u_x,
$$
we obtain
\begin{equation}
\begin{array}{ll}
\di \Big[\f{\r(u-\bar u)^2}{2}+\r\Psi(\r,\bar\r)\Big]_t+H_{1x}(t,x)
+\mu_\v(\r)\Big[(u-\bar u)_x\Big]^2+\bar
u_x\Big[\r^\g-\bar\r^\g-\g\bar\r^{\g-1}(\r-\bar\r)\Big]\\\di
+\r(u-\bar u)^2\bar u_x=\Big[\mu_\v(\r)\Big]_x\bar u_x(u-\bar
u)+\mu_\v(\r)\bar u_{xx}(u-\bar u).
\end{array}\label{(3.9)}
\end{equation}

{\it Step 2. Entropy Equality}

Rewrite \eqref{(3.5)} as
\begin{equation}
\r(u-\bar u)_t+\r u(u-\bar u)_x+(p(\r)-p(\bar\r))_x+(\r-\bar\r)\bar
u_t+(\r u-\bar\r\bar u)\bar u_x=\Big[(\r^{\a-1}+\v\r^{\t-1})\r
u_x\Big]_x.\label{(3.10)}
\end{equation}
Note that
\begin{equation}
\di \Big[(\r^{\a-1}+\v\r^{\t-1})\r u_x\Big]_x \di
=-\r(\varphi_\v^{\a,\t}(\r))_{xt}-\r
u(\varphi_\v^{\a,\t}(\r))_{xx},\label{(3.11)}
\end{equation}
where $\varphi_\v^{\a,\t}(\r)$ with $\t=\f12$ is defined by
$$
\varphi_\v^{\a,\t}(\r)=\left\{
\begin{array}{ll}
\di \f{\r^{\a-1}}{\a-1}+\v\f{\r^{\t-1}}{\t-1},&\di {\rm
if}~~\a\neq1,\a>0,\\
\di \ln\r+\v\f{\r^{\t-1}}{\t-1},&\di {\rm if}~~\a=1.
\end{array}
\right.
$$
Thus \eqref{(3.10)} becomes
\begin{equation}
\r(u-\bar u)_t+\r u(u-\bar u)_x+(p(\r)-p(\bar\r))_x+(\r-\bar\r)\bar
u_t+(\r u-\bar\r\bar u)\bar u_x=-\r(\varphi_\v^{\a,\t}(\r))_{xt}-\r
u(\varphi_\v^{\a,\v}(\r))_{xx}.\label{(3.12)}
\end{equation}
Multiplying \eqref{(3.12)} by $(\varphi_\v^{\a,\t}(\r))_x$ shows
that
\begin{equation}
\begin{array}{ll}
\di \Big[\f{\r(\varphi_\v^{\a,\t}(\r))_x^2}{2}\Big]_t+ \Big[\f{\r
u(\varphi_\v^{\a,\t}(\r))_x^2}{2}\Big]_x+\Big[\r(u-\bar
u)(\varphi_\v^{\a,\t}(\r))_x\Big]_t +\Big[\r u(u-\bar
u)(\varphi_\v^{\a,\t}(\r))_x\Big]_x\\
\di -(u-\bar u)\Big[\r(\varphi_\v^{\a,\t}(\r))_{xt}+\r
u(\varphi_\v^{\a,\t}(\r))_{xx}\Big]+(\varphi_\v^{\a,\t}(\r))_x(p(\r)-p(\bar\r))_x\\
\di +(\varphi_\v^{\a,\t}(\r))_x\Big[(\r-\bar\r)\bar u_t+(\r
u-\bar\r\bar u)\bar u_x\Big]=0.
\end{array}\label{(3.13)}
\end{equation}
Combining \eqref{(3.12)} with \eqref{(3.13)} yields
\begin{equation}
\begin{array}{ll}
\di \Big\{\f12\r\left[(u-\bar
u)+(\varphi_\v^{\a,\t}(\r))_x\right]^2\Big\}_t+ \Big\{\f12\r
u\left[(u-\bar
u)+(\varphi_\v^{\a,\t}(\r))_x\right]^2\Big\}_x+(u-\bar u)(p(\r)-p(\bar\r))_x\\
\di +(\varphi_\v^{\a,\t}(\r))_x(p(\r)-p(\bar\r))_x+(u-\bar
u)\Big[(\r-\bar\r)\bar u_t+(\r
u-\bar\r\bar u)\bar u_x\Big]\\
\di +(\varphi_\v^{\a,\t}(\r))_x\Big[(\r-\bar\r)\bar u_t+(\r
u-\bar\r\bar u)\bar u_x\Big]=0.
\end{array}\label{(3.14)}
\end{equation}

{\it Step 3. A Priori Estimates}

It follows from \eqref{(3.7)} and \eqref{(3.14)} that
\begin{equation}
\begin{array}{ll}
\di \Big\{\f12\r\left[(u-\bar
u)+(\varphi_\v^{\a,\t}(\r))_x\right]^2+\r\Psi(\r,\bar\r)\Big\}_t+
\Big\{\f12\r u\left[(u-\bar
u)+(\varphi_\v^{\a,\t}(\r))_x\right]^2+\r u\Psi(\r,\bar\r)\\
\di +(u-\bar u)(p(\r)-p(\bar\r))\Big\}_x+\bar
u_x\Big[p(\r)-p(\bar\r)-p^\prime(\bar\r)(\r-\bar\r)\Big]+\r(u-\bar u)^2\bar u_x\\
\di +(\varphi_\v^{\a,\t}(\r))_x\Big[(\r-\bar\r)\bar u_t+(\r
u-\bar\r\bar u)\bar u_x+p(\r)_x-p(\bar\r)_x\Big]=0.
\end{array}\label{(3.15)}
\end{equation}
Now we deal with the last term on the left hand side of
\eqref{(3.15)}. Note that
\begin{equation}
\di (\r-\bar\r)\bar u_t+(\r u-\bar\r\bar u)\bar
u_x+p(\r)_x-p(\bar\r)_x \di =\r(u-\bar u)\bar u_x+\Big[p(\r)_x-\f{\r
p(\bar\r)_x}{\bar\r}\Big],\label{(3.16)}
\end{equation}
and
\begin{equation}
(\varphi_\v^{\a,\t}(\r))_x=\r^{\a-2}\r_x+\v\r^{\t-2}\r_x.\label{(3.17)}
\end{equation}
Thus
\begin{equation}
\begin{array}{ll}
&\di (\varphi_\v^{\a,\t}(\r))_x\Big[(\r-\bar\r)\bar u_t+(\r
u-\bar\r\bar u)\bar u_x+p(\r)_x-p(\bar\r)_x\Big]\\
&\di =\Big(\f{\r^\a}{\a}+\v\f{\r^\t}{\t}\Big)_x(u-\bar u)\bar
u_x+(\r^{\a-2}\r_x+\v\r^{\t-2}\r_x)\Big[p(\r)_x-\f{\r
p(\bar\r)_x}{\bar\r}\Big].\label{(3.18)}
\end{array}
\end{equation}
Direct computations show
\begin{equation}
\begin{array}{ll}
&\di \r^{\a-2}\r_x\Big[p(\r)_x-\f{\r p(\bar\r)_x}{\bar\r}\Big]\\
&\di
=\f{4\g}{(\a+\g-1)^2}\Big[(\r^{\f{\a+\g-1}{2}}-\bar\r^{\f{\a+\g-1}{2}})_x\Big]^2+\Big[\f{8\g}{(\a+\g-1)^2}(\bar\r^{\f{\a+\g-1}{2}})_x(\r^{\f{\a+\g-1}{2}}-\bar\r^{\f{\a+\g-1}{2}})\\
&\di \quad-\f{2\g}{\a(\a+\g-1)}(\bar\r^{\f{\a+\g-1}{2}})_x\bar\r^{\f{\g-\a-1}{2}}(\r^\a-\bar\r^\a)\Big]_x-\f{8\g}{(\a+\g-1)^2}(\bar\r^{\f{\a+\g-1}{2}})_{xx}(\r^{\f{\a+\g-1}{2}}-\bar\r^{\f{\a+\g-1}{2}})\\
&\di \quad+\f{2\g}{\a(\a+\g-1)}\Big[(\bar\r^{\f{\a+\g-1}{2}})_x\bar\r^{\f{\g-\a-1}{2}}\Big]_x(\r^\a-\bar\r^\a).\\
\end{array}\label{(3.19)}
\end{equation}
and
\begin{equation}
\begin{array}{ll}
&\di \r^{\t-2}\r_x\Big[p(\r)_x-\f{\r p(\bar\r)_x}{\bar\r}\Big]\\
&\di
=\f{4\g}{(\t+\g-1)^2}\Big[(\r^{\f{\t+\g-1}{2}}-\bar\r^{\f{\t+\g-1}{2}})_x\Big]^2+\Big[\f{8\g}{(\t+\g-1)^2}(\bar\r^{\f{\t+\g-1}{2}})_x(\r^{\f{\t+\g-1}{2}}-\bar\r^{\f{\t+\g-1}{2}})\\
&\di \quad-\f{2\g}{\t(\t+\g-1)}(\bar\r^{\f{\t+\g-1}{2}})_x\bar\r^{\f{\g-\t-1}{2}}(\r^\t-\bar\r^\t)\Big]_x-\f{8\g}{(\t+\g-1)^2}(\bar\r^{\f{\t+\g-1}{2}})_{xx}(\r^{\f{\t+\g-1}{2}}-\bar\r^{\f{\t+\g-1}{2}})\\
&\di \quad+\f{2\g}{\t(\t+\g-1)}\Big[(\bar\r^{\f{\t+\g-1}{2}})_x\bar\r^{\f{\g-\t-1}{2}}\Big]_x(\r^\t-\bar\r^\t).\\
\end{array}\label{(3.20)}
\end{equation}
Substituting \eqref{(3.18)}-\eqref{(3.20)} into \eqref{(3.15)} gives
\begin{equation}
\begin{array}{ll}
\di \Big\{\f12\r\left[(u-\bar
u)+(\varphi_\v^{\a,\t}(\r))_x\right]^2+\r\Psi(\r,\bar\r)\Big\}_t+
H_{2x}(t,x)+\bar
u_x\Big[p(\r)-p(\bar\r)-p^\prime(\bar\r)(\r-\bar\r)\Big]\\
\di +\r(u-\bar u)^2\bar
u_x+\Big(\f{\r^\a}{\a}+\v\f{\r^\t}{\t}\Big)_x(u-\bar u)\bar
u_x+\f{4\g}{(\a+\g-1)^2}\Big[(\r^{\f{\a+\g-1}{2}}-\bar\r^{\f{\a+\g-1}{2}})_x\Big]^2\\
\di
+\v\f{4\g}{(\t+\g-1)^2}\Big[(\r^{\f{\t+\g-1}{2}}-\bar\r^{\f{\t+\g-1}{2}})_x\Big]^2=\f{8\g}{(\a+\g-1)^2}(\bar\r^{\f{\a+\g-1}{2}})_{xx}(\r^{\f{\a+\g-1}{2}}-\bar\r^{\f{\a+\g-1}{2}})\\
\di
+\v\f{8\g}{(\t+\g-1)^2}(\bar\r^{\f{\t+\g-1}{2}})_{xx}(\r^{\f{\t+\g-1}{2}}-\bar\r^{\f{\t+\g-1}{2}})-\f{2\g}{\a(\a+\g-1)}\Big[(\bar\r^{\f{\a+\g-1}{2}})_x\bar\r^{\f{\g-\a-1}{2}}\Big]_x(\r^\a-\bar\r^\a)\\
\di
-\v\f{2\g}{\t(\t+\g-1)}\Big[(\bar\r^{\f{\t+\g-1}{2}})_x\bar\r^{\f{\g-\t-1}{2}}\Big]_x(\r^\t-\bar\r^\t),
\end{array}\label{(3.21)}
\end{equation}
where
\begin{equation}
\begin{array}{ll}
\di H_2(t,x)=&\di\f12\r u\left[(u-\bar
u)+(\varphi_\v^{\a,\t}(\r))_x\right]^2+\r u\Psi(\r,\bar\r)+(u-\bar
u)(p(\r)-p(\bar\r))\\
&\di
+\f{8\g}{(\a+\g-1)^2}(\bar\r^{\f{\a+\g-1}{2}})_x(\r^{\f{\a+\g-1}{2}}-\bar\r^{\f{\a+\g-1}{2}})\\
&\di
-\f{2\g}{\a(\a+\g-1)}(\bar\r^{\f{\a+\g-1}{2}})_x\bar\r^{\f{\g-\a-1}{2}}(\r^\a-\bar\r^\a)\\
&\di +\v\f{8\g}{(\t+\g-1)^2}(\bar\r^{\f{\t+\g-1}{2}})_x(\r^{\f{\t+\g-1}{2}}-\bar\r^{\f{\t+\g-1}{2}})\\
&\di
-\v\f{2\g}{\t(\t+\g-1)}(\bar\r^{\f{\t+\g-1}{2}})_x\bar\r^{\f{\g-\t-1}{2}}(\r^\t-\bar\r^\t).
\end{array}\label{(3.22)}
\end{equation}
 Multiplying \eqref{(3.21)} by $\a$ and  then
adding up to \eqref{(3.9)} and noticing that
$\Big[\mu_\v(\r)\Big]_x=(\r^\a)_x+\v(\r^\t)_x$ in the right hand
side of \eqref{(3.9)}, one can get
\begin{equation}
\begin{array}{ll}
\di \Big\{\f{\a}{2}\r\left[(u-\bar
u)+(\varphi_\v^{\a,\t}(\r))_x\right]^2+\f{\r(u-\bar
u)^2}{2}+(\a+1)\r\Psi(\r,\bar\r)\Big\}_t+ \Big[\a
H_{2}(t,x)+H_1(t,x)\Big]_x\\
\di +(\a+1)\bar
u_x\Big[p(\r)-p(\bar\r)-p^\prime(\bar\r)(\r-\bar\r)\Big]+(\a+1)\r(u-\bar
u)^2\bar u_x+(\r^\a+\v\r^\t)\Big[(u-\bar u)_x\Big]^2\\
\di
+\f{4\a\g}{(\a+\g-1)^2}\Big[(\r^{\f{\a+\g-1}{2}}-\bar\r^{\f{\a+\g-1}{2}})_x\Big]^2
+\v\f{4\a\g}{(\t+\g-1)^2}\Big[(\r^{\f{\t+\g-1}{2}}-\bar\r^{\f{\t+\g-1}{2}})_x\Big]^2\\
\di=\r^\a\bar u_{xx}(u-\bar u)+\v\Big[\r^\t\bar u_{xx}(u-\bar u)+(1-\f\a\t)(\r^\t)_x(u-\bar u)\bar u_x\Big]\\
\di
+\f{8\a\g}{(\a+\g-1)^2}(\bar\r^{\f{\a+\g-1}{2}})_{xx}(\r^{\f{\a+\g-1}{2}}-\bar\r^{\f{\a+\g-1}{2}})
+\v\f{8\a\g}{(\t+\g-1)^2}(\bar\r^{\f{\t+\g-1}{2}})_{xx}(\r^{\f{\t+\g-1}{2}}-\bar\r^{\f{\t+\g-1}{2}})
\\
\di-\f{2\g}{(\a+\g-1)}\Big[(\bar\r^{\f{\a+\g-1}{2}})_x\bar\r^{\f{\g-\a-1}{2}}\Big]_x(\r^\a-\bar\r^\a)
-\v\f{2\a\g}{\t(\t+\g-1)}\Big[(\bar\r^{\f{\t+\g-1}{2}})_x\bar\r^{\f{\g-\t-1}{2}}\Big]_x(\r^\t-\bar\r^\t).
\end{array}\label{(3.22)}
\end{equation}
Integrating \eqref{(3.22)} over $[0,t]\times\mathbf{R}$ with respect
to $t,x$ gives
\begin{equation}
\begin{array}{ll}
\di \int_{\mathbf{R}}\Big\{\f{\a}{2}\r\left[(u-\bar
u)+(\varphi_\v^{\a,\t}(\r))_x\right]^2+\f{\r(u-\bar
u)^2}{2}+(\a+1)\r\Psi(\r,\bar\r)\Big\}(t,x)dx \\
 \di +\int_0^t\int_{\mathbf{R}}\Big\{(\a+1)\bar
u_x\Big[p(\r)-p(\bar\r)-p^\prime(\bar\r)(\r-\bar\r)\Big]+(\a+1)\r(u-\bar
u)^2\bar u_x\\
\di+(\r^\a+\v\r^\t)\Big[(u-\bar u)_x\Big]^2
+\f{4\a\g}{(\a+\g-1)^2}\Big[(\r^{\f{\a+\g-1}{2}}-\bar\r^{\f{\a+\g-1}{2}})_x\Big]^2\\
\di
+\v\f{4\a\g}{(\t+\g-1)^2}\Big[(\r^{\f{\t+\g-1}{2}}-\bar\r^{\f{\t+\g-1}{2}})_x\Big]^2\Big\}dxd\tau\\
\di = \int_{\mathbf{R}}\Big\{\f{\a}{2}\r_0\left[(u_0-\bar
u_0)+(\varphi_\v^{\a,\t}(\r_0))_x\right]^2+\f{\r_0(u_0-\bar
u_0)^2}{2}+(\a+1)\r_0\Psi(\r_0,\bar\r_0)\Big\}dx+I,\\
\end{array}
\end{equation}
where
\begin{equation}
\begin{array}{ll}
\di I=\int_0^t\int_{\mathbf{R}}\Big\{\r^\a\bar u_{xx}(u-\bar u)+\v\Big[\r^\t\bar u_{xx}(u-\bar u)+(1-\f\a\t)(\r^\t)_x(u-\bar u)\bar u_x\Big]\\
\di
+\f{8\a\g}{(\a+\g-1)^2}(\bar\r^{\f{\a+\g-1}{2}})_{xx}(\r^{\f{\a+\g-1}{2}}-\bar\r^{\f{\a+\g-1}{2}})
+\v\f{8\a\g}{(\t+\g-1)^2}(\bar\r^{\f{\t+\g-1}{2}})_{xx}(\r^{\f{\t+\g-1}{2}}-\bar\r^{\f{\t+\g-1}{2}})
\\
\di-\f{2\g}{(\a+\g-1)}\Big[(\bar\r^{\f{\a+\g-1}{2}})_x\bar\r^{\f{\g-\a-1}{2}}\Big]_x(\r^\a-\bar\r^\a)
-\v\f{2\a\g}{\t(\t+\g-1)}\Big[(\bar\r^{\f{\t+\g-1}{2}})_x\bar\r^{\f{\g-\t-1}{2}}\Big]_x(\r^\t-\bar\r^\t)\Big\}dxd\tau\\
\di :=\sum_{i=1}^6I_i.
\end{array}\label{I}
\end{equation}
We now estimate the right hand side of \eqref{I} terms by terms.
First,
\begin{equation}
\begin{array}{ll}
\di I_1&\di =\int_0^t\int_{\mathbf{R}}\r^\a\bar
u_{xx}(u-\bar u)dxd\tau\\
&\di =\int_0^t\int_{\mathbf{R}}\sqrt{\r}(u-\bar u)\r^{\a-\f12}\bar
u_{xx}dxd\tau\\
&\di =\int_0^t\int_{\mathbf{R}}\sqrt{\r}(u-\bar u)\r^{\a-\f12}\bar
u_{xx}[{\bf 1}|_{\{0\leq\r\leq2\r_+\}}+{\bf
1}|_{\{\r\geq2\r_+\}}]dxd\tau\\
&\di :=I_{11}+I_{12},
\end{array}\label{I1}
\end{equation}
where and in the sequel ${\bf 1}|_\Omega$ denotes the characteristic
function of a set $\Omega\subset(0,t)\times\mathbf{R}$.

Rewrite $I_{12}$ as
\begin{equation}
\begin{array}{ll}
\di I_{12}=\int_0^t\int_{\mathbf{R}}\sqrt{\r}(u-\bar u)\bar
u_{xx}[(\r^{\a-\f12}-\bar\r^{\a-\f12})+\bar\r^{\a-\f12}]{\bf
1}|_{\{\r\geq2\r_+\}}dxd\tau\\
\di\quad:=I_{12}^1+I_{12}^2.
\end{array}\label{I12}
\end{equation}
Using Lemma \ref{rarefaction-1} (and its Remark 2.1) and noting that
$\a\geq 1$,  one has
\begin{equation}
\begin{array}{ll}
\di I_{11}+I_{12}^2&\di \leq C\int_0^t\|\sqrt\r(u-\bar
u)\|_{L^2(\mathbf{R})}\|\bar
u_{xx}\|_{L^2(\mathbf{R})}d\tau\\
&\di\leq C\sup_{t\in[0,T]}\|\sqrt\r(u-\bar
u)\|_{L^2(\mathbf{R})}\int_0^t\|\bar
u_{xx}\|_{L^2(\mathbf{R})}d\tau\\
&\di \leq C\sup_{t\in[0,T]}\|\sqrt\r(u-\bar u)\|_{L^2(\mathbf{R})},
\end{array}\label{I11}
\end{equation}
and
\begin{equation}
\begin{array}{ll}
\di I_{12}^1\leq C\sup_{t\in[0,T]}\|\sqrt{\r}(u-\bar
u)\|_{L^2(\mathbf{R})}\sup_{t\in[0,T]}\|(\r^{\a-\f12}-\bar\r^{\a-\f12}){\bf
1}|_{\{\r\geq2\r_+\}}\|_{L^2(\mathbf{R})}\cdot\int_0^t\|\bar u_{xx}\|_{L^\i(\mathbf{R})}d\tau\\
\di\quad\leq C\eta^{\f{2}{4q+1}}\sup_{t\in[0,T]}\|\sqrt{\r}(u-\bar
u)\|_{L^2(\mathbf{R})}\sup_{t\in[0,T]}\|(\r^{\a-\f12}-\bar\r^{\a-\f12}){\bf 1}|_{\{\r\geq2\r_+\}}\|_{L^2(\mathbf{R})}\\
\di\quad\leq
C\eta^{\f{2}{4q+1}}\Big[\sup_{t\in[0,T]}\|\sqrt\r(u-\bar
u)\|^2_{L^2(\mathbf{R})}+\sup_{t\in[0,T]}\|(\r^{\a-\f12}-\bar\r^{\a-\f12}){\bf
1}|_{\{\r\geq2\r_+\}}\|^2_{L^2(\mathbf{R})}\Big].
\end{array}\label{I121}
\end{equation}
Note that if $\a$ and $\g$ satisfy
\begin{equation}
1\leq\a\leq \f{\g+1}{2},
\end{equation}
then $2(\a-\f12)\leq\g$, and then
\begin{equation}
\begin{array}{ll}
\di
\lim_{\r\rightarrow+\i}\f{(\r^{\a-\f12}-\bar\r^{\a-\f12})^2}{\r\Psi(\r,\bar\r)}&\di
=\lim_{\r\rightarrow+\i}\f{(\g-1)(\r^{\a-\f12}-\bar\r^{\a-\f12})^2}{\r^\g-\bar\r^\g-\g\bar\r^{\g-1}(\r-\bar\r)}\leq
C.
\end{array}\label{I1+}
\end{equation}
Thus if $1 \leq\a\leq \f{\g+1}{2}$, then
\begin{equation}
\sup_{t\in[0,T]}\|(\r^{\a-\f12}-\bar\r^{\a-\f12}){\bf
1}|_{\{\r\geq2\r_+\}}\|^2_{L^2(\mathbf{R})}\leq
C\sup_{t\in[0,T]}\|\r\Psi(\r,\bar\r)\|_{L^1(\mathbf{R})},\label{I1++}
\end{equation}
for some uniform constant $C>0$.

Substituting \eqref{I12}, \eqref{I11}, \eqref{I121} and \eqref{I1++}
into \eqref{I1} yields
\begin{equation}
I_1\leq C\eta^{\f{2}{4q+1}}\Big[\sup_{t\in[0,T]}\|\sqrt\r(u-\bar
u)\|^2_{L^2(\mathbf{R})}+\sup_{t\in[0,T]}\|\r\Psi(\r,\bar\r)\|_{L^1(\mathbf{R})}\Big]+C_\eta.\label{I1E}
\end{equation}
Next, $I_2$ can be rewritten as
\begin{equation}
\begin{array}{ll}
I_2&\di =\v\int_0^t\int_{\mathbf{R}}\Big[\r^\t\bar
u_{xx}(u-\bar u)+(1-\f\a\t)(\r^\t)_x(u-\bar u)\bar u_x\Big]dxd\tau\\
&\di =\v\int_0^t\int_{\mathbf{R}}\Big[\f\a\t\r^\t\bar u_{xx}(u-\bar
u)-(1-\f\a\t)\r^\t(u-\bar u)_x\bar u_x\Big]dxd\tau\\
&\di :=I_{21}+I_{22},
\end{array}
\label{I2}
\end{equation}
 First, since $\t=\f12$, it follows that
 \begin{equation}
\begin{array}{ll}
\di I_{21}&\di=\v\int_0^t\int_{\mathbf{R}}2\a\bar
u_{xx}\sqrt\r(u-\bar u)dxd\tau\\
&\di\leq C\v\sup_{t\in[0,T]}\|\sqrt{\r}(u-\bar
u)\|_{L^2(\mathbf{R})}\int_0^t\|\bar u_{xx}\|_{L^2(\mathbf{R})}d\tau\\
&\di\leq C\v\sup_{t\in[0,T]}\|\sqrt{\r}(u-\bar
u)\|_{L^2(\mathbf{R})}.
\end{array}\label{I21}
\end{equation}
On the other hand,
\begin{equation}
\begin{array}{ll}
\di I_{22}\leq \f\v4 \int_0^t\int_{\mathbf{R}}\r^\t[(u-\bar
u)_x]^2dxd\tau+C\v\int_0^t\int_{\mathbf{R}}\r^\t\bar u_x^2dxd\tau,
\end{array}\label{I22}
\end{equation}
while
\begin{equation}
\begin{array}{ll}
\di \v\int_0^t\int_{\mathbf{R}}\r^\t\bar u_x^2dxd\tau&\di
=\v\int_0^t\int_{\mathbf{R}}\r^\t[{\bf
1}|_{\{0\leq\r\leq2\r_+\}}+{\bf 1}|_{\{\r\geq2\r_+\}}]\bar
u_x^2dxd\tau\\
&\di =\v\int_0^t\int_{\mathbf{R}}\Big\{\r^\t{\bf
1}|_{\{0\leq\r\leq2\r_+\}}+\big[(\r^\t-\bar\r^\t)+\bar\r^\t\big]{\bf
1}|_{\{\r\geq2\r_+\}}\Big\}\bar
u_x^2dxd\tau\\
&\di \leq C\v\ln(1+T)+C\v\sup_{t\in[0,T]}\|(\r^\t-\bar\r^\t){\bf
1}|_{\{\r\geq2\r_+\}}\|_{L^1(\mathbf{R})}\int_0^t\|\bar
u_x\|_{L^\i(\mathbf{R})}^2d\tau\\
&\di \leq
C\v\ln(1+T)+C\v\sup_{t\in[0,T]}\|\r\Psi(\r,\bar\r)\|_{L^1(\mathbf{R})},
\end{array}\label{I221}
\end{equation}
where in the last inequality we have used the fact that
$$
\lim_{\r\rightarrow+\i}\f{(\r^\t-\bar\r^\t){\bf
1}|_{\{\r\geq2\r_+\}}}{\r\Psi(\r,\bar\r)}=0,
$$
since $\t=\f12<1< \g.$

Substituting the estimations \eqref{I21}-\eqref{I221} into
\eqref{I2}, one can get
\begin{equation}
\begin{array}{ll}
I_2  &\di\leq  C\v\ln(1+T)+\f\v4
\int_0^t\int_{\mathbf{R}}\r^\t[(u-\bar
u)_x]^2dxd\tau+C\v\sup_{t\in[0,T]}\|\r\Psi(\r,\bar\r)\|_{L^1(\mathbf{R})}\\
&\di +C\v\sup_{t\in[0,T]}\|\sqrt{\r}(u-\bar u)\|_{L^2(\mathbf{R})}.
\end{array}  \label{I2E}
\end{equation}
It follows from the fact that
\begin{equation}
\lim_{\r\rightarrow
0}\f{|\r^{\f{\a+\g-1}{2}}-\bar\r^{\f{\a+\g-1}{2}}|^{\f{2\g}{\a+\g-1}}}{\r\Psi(\r,\bar\r)}=1,
\end{equation}
that for any $\epsilon>0$, there exists $\d_\epsilon>0,$ such that
if $0\leq \r\leq\d_\epsilon$, then
\begin{equation}
|\f{|\r^{\f{\a+\g-1}{2}}-\bar\r^{\f{\a+\g-1}{2}}|^{\f{2\g}{\a+\g-1}}}{\r\Psi(\r,\bar\r)}-1|<\epsilon.\label{fa1}
\end{equation}
Fix $\epsilon=\f12,$ then there exists $\d_{\f12}>0,$ such that if
$0\leq \r\leq\d_{\f12}$, then
\begin{equation}
|\f{|\r^{\f{\a+\g-1}{2}}-\bar\r^{\f{\a+\g-1}{2}}|^{\f{2\g}{\a+\g-1}}}{\r\Psi(\r,\bar\r)}-1|<\f12,\label{fa2}
\end{equation}
thus for any $\bar\r\geq0$,
\begin{equation}
\f12\r\Psi(\r,\bar\r)\leq
|\r^{\f{\a+\g-1}{2}}-\bar\r^{\f{\a+\g-1}{2}}|^{\f{2\g}{\a+\g-1}}\leq
\f32\r\Psi(\r,\bar\r),\quad {\rm if} ~~0\leq
\r\leq\d_{\f12}.\label{fa3}
\end{equation}
Similarly, it follows from the fact that
\begin{equation}
\lim_{\bar\r\rightarrow
0}\f{|\r^{\f{\a+\g-1}{2}}-\bar\r^{\f{\a+\g-1}{2}}|^{\f{2\g}{\a+\g-1}}}{\r\Psi(\r,\bar\r)}=1,
\end{equation} that there exists $\bar\d_{\f12}>0,$ such that if
$0\leq \bar\r\leq\bar\d_{\f12}$, then
\begin{equation}
|\f{|\r^{\f{\a+\g-1}{2}}-\bar\r^{\f{\a+\g-1}{2}}|^{\f{2\g}{\a+\g-1}}}{\r\Psi(\r,\bar\r)}-1|<\f12,\label{fa4}
\end{equation}
thus one can choose $\nu< \bar\d_{\f12}$ such that for any
$\r\geq0,$
\begin{equation}
\f12\r\Psi(\r,\bar\r)\leq
|\r^{\f{\a+\g-1}{2}}-\bar\r^{\f{\a+\g-1}{2}}|^{\f{2\g}{\a+\g-1}}\leq
\f32\r\Psi(\r,\bar\r),\quad {\rm if} ~~\nu\leq
\bar\r\leq\bar\d_{\f12}.\label{fa5}
\end{equation}
The term $I_3$ can be estimated as follows. Since
$$
\begin{array}{ll}
\di (\bar\r^{\f{\a+\g-1}{2}})_{xx}&\di
=(\f{\a+\g-1}{2}\bar\r^{\f{\a+\g-3}{2}}\bar\r_x)_x=\f{\a+\g-1}{2\sqrt\g}(\bar\r^{\f\a2}\bar
u_x)_x\\
&\di =\f{\a+\g-1}{2\sqrt\g}\bar\r^{\f\a2}\bar
u_{xx}+\f{\a(\a+\g-1)}{4\g}\bar\r^{\f{\a+1-\g}{2}}\bar u_x^2,
\end{array}
$$
and
$$
\a+1-\g\geq 2-\g\geq0,
$$
one can rewrite $I_3$ as
\begin{equation}
\begin{array}{ll}
\di
I_3&\di=\int_0^t\int_{\mathbf{R}}\Big(\f{4\a\sqrt\g}{\a+\g-1}\bar\r^{\f\a2}\bar
u_{xx}
+\f{8\a^2}{\a+\g-1}\bar
u_x^2\bar\r^{\f{\a+1-\g}{2}}\Big)(\r^{\f{\a+\g-1}{2}}-\bar\r^{\f{\a+\g-1}{2}})dxd\tau\\[3mm]
&\di \leq C\int_0^t\int_{\mathbf{R}}|(\bar u_{xx},\bar
u_x^2)||(\r^{\f{\a+\g-1}{2}}-\bar\r^{\f{\a+\g-1}{2}})|dxd\tau\\
&\di =C\int_0^t\int_{\mathbf{R}}|(\bar u_{xx},\bar
u_x^2)||(\r^{\f{\a+\g-1}{2}}-\bar\r^{\f{\a+\g-1}{2}})|\big({\bf
1}|_{\{0\leq\r\leq\d_{\f12}\}}+{\bf
1}|_{\{\d_{\f12}\leq\r\leq2\r_+,~\nu\leq
\bar\r\leq\bar\d_{\f12}\}}\\
&\qquad\qquad\qquad\di+{\bf
1}|_{\{\d_{\f12}\leq\r\leq2\r_+,~\bar\d_{\f12}\leq
\bar\r\leq\r_+\}}+{\bf
1}|_{\{\r\geq2\r_+\}}\big)dxd\tau\\
 &\di:=I_{31}+I_{32}+I_{33}+I_{34}.
\end{array}
\end{equation}
Direct computations lead to
\begin{equation}
\begin{array}{ll}
\di I_{31}&\di \leq C\int_0^t\int_{\mathbf{R}}\|(\bar u_{xx},\bar
u_x^2)\|_{L^{\f{2\g}{\g+1-\a}}(\mathbf{R})}\|(\r^{\f{\a+\g-1}{2}}-\bar\r^{\f{\a+\g-1}{2}}){\bf
1}|_{\{0\leq\r\leq\d_{\f12}\}}\|_{L^{\f{2\g}{\a+\g-1}}(\mathbf{R})}dxd\tau\\[3mm]
&\di\le
C\sup_{t\in[0,T]}\|(\r^{\f{\a+\g-1}{2}}-\bar\r^{\f{\a+\g-1}{2}}){\bf
1}|_{\{0\leq\r\leq\d_{\f12}\}}\|_{L^{\f{2\g}{\a+\g-1}}(\mb{R})}\int_0^t\|(\bar
u_{xx},\bar
u_x^2)\|_{L^{\f{2\g}{\g+1-\a}}(\mathbf{R})}d\tau\\[3mm]
&\di \stackrel{\eqref{fa3}}{\leq}
C\sup_{t\in[0,T]}\|\r\Psi(\r,\bar\r)\|^{\f{\a+\g-1}{2\g}}_{L^1(\mathbf{R})}\\
&\di \leq
\f{\a+1}{8}\sup_{t\in[0,T]}\|\r\Psi(\r,\bar\r)\|_{L^1(\mathbf{R})}+C_\a.
\end{array}
\end{equation}
Similarly, due to \eqref{fa5}, one has
\begin{equation}
I_{32}\leq
\f{\a+1}{8}\sup_{t\in[0,T]}\|\r\Psi(\r,\bar\r)\|_{L^1(\mathbf{R})}+C_\a.
\end{equation}
On the other hand,
\begin{equation}
\begin{array}{ll}
I_{33}&\di\leq C\int_0^t\int_{\mathbf{R}}\|(\bar u_{xx},\bar
u_x^2)\|_{L^2(\mathbf{R})}\|(\r^{\f{\a+\g-1}{2}}-\bar\r^{\f{\a+\g-1}{2}}){\bf
1}|_{\{\d_{\f12}\leq\r\leq2\r_+,~\bar\d_{\f12}\leq
\bar\r\leq\r_+\}}\|_{L^2(\mathbf{R})}dxd\tau\\[3mm]
&\di\le
C\sup_{t\in[0,T]}\|(\r^{\f{\a+\g-1}{2}}-\bar\r^{\f{\a+\g-1}{2}}){\bf
1}|_{\{\d_{\f12}\leq\r\leq2\r_+,~\bar\d_{\f12}\leq
\bar\r\leq\r_+\}}\|_{L^2(\mb{R})}\int_0^t\|(\bar u_{xx},\bar
u_x^2)\|_{L^2(\mathbf{R})}d\tau\\[3mm]
&\di \leq
C\sup_{t\in[0,T]}\|\r\Psi(\r,\bar\r)\|^{\f12}_{L^1(\mathbf{R})},
\end{array}
\end{equation}
where one has used the fact that
$$
\di \f{(\r^{\f{\a+\g-1}{2}}-\bar\r^{\f{\a+\g-1}{2}})^2{\bf
1}|_{\{\d_{\f12}\leq\r\leq2\r_+,~\bar\d_{\f12}\leq
\bar\r\leq\r_+\}}}{\r\Psi(\r,\bar\r)}\leq C.
$$
Moreover, $I_{34}$ can be estimated as
\begin{equation}
\begin{array}{ll}
\di I_{34}&\di\leq C\int_0^t\int_{\mathbf{R}}\|(\bar u_{xx},\bar
u_x^2)\|_{L^\i(\mathbf{R})}\|(\r^{\f{\a+\g-1}{2}}-\bar\r^{\f{\a+\g-1}{2}}){\bf
1}_{\{\r\geq2\r_+\}}\|_{L^1(\mathbf{R})}dxd\tau\\[3mm]
&\di\le
C\sup_{t\in[0,T]}\|(\r^{\f{\a+\g-1}{2}}-\bar\r^{\f{\a+\g-1}{2}}){\bf
1}_{\{\r\geq2\r_+\}}\|_{L^1(\mb{R})} \int_0^t\|(\bar u_{xx},\bar
u_x^2)\|_{L^\i(\mathbf{R})}d\tau\\[3mm]
&\di \leq
C\eta^\f{2}{4q+1}\sup_{t\in[0,T]}\|\r\Psi(\r,\bar\r)\|_{L^1(\mathbf{R})},
\end{array}
\end{equation}
due to the facts that
$$
\lim_{\r\rightarrow+\i}\f{|\r^{\f{\a+\g-1}{2}}-\bar\r^{\f{\a+\g-1}{2}}|{\bf
1}_{\{\r\geq2\r_+\}}}{\r\Psi(\r,\bar\r)}\leq C,
$$
since $\a\leq \f{\g+1}{2}$ implies
$$
\f{\a+\g-1}{2}\leq \g, \quad{\rm i.\ e.},\quad \a\leq \g+1.
$$
Now we turn to the term $I_5$. First,
$$
\begin{array}{ll}
\di
\Big[(\bar\r^{\f{\a+\g-1}{2}})_x\bar\r^{\f{\g-\a-1}{2}}\Big]_x&\di
=\f{\a+\g-1}{2\sqrt{\g}}(\bar\r^{\f{\g-1}{2}}\bar
u_x)_x\\
&\di=\f{\a+\g-1}{2\sqrt{\g}}\bar\r^{\f{\g-1}{2}}\bar
u_{xx}+\f{(\a+\g-1)(\g-1)}{4\g}(\bar u_x)^2.
\end{array}
$$
Thus
$$
\begin{array}{ll}
\di I_5&\di
=\int_0^t\int_{\mathbf{R}}\f{2\g}{(\a+\g-1)}\Big[(\bar\r^{\f{\a+\g-1}{2}})_x\bar\r^{\f{\g-\a-1}{2}}\Big]_x(\r^\a-\bar\r^\a)dxd\tau\\[3mm]
 &\di \leq C\int_0^t\int_{\mathbf{R}}|(\bar u_{xx},\bar u_x^2)| |(\r^\a-\bar\r^\a)|({\bf
1}|_{\{0\leq\r\leq2\r_+\}}+{\bf
1}|_{\{\r\geq2\r_+\}})dxd\tau\\
&\di :=I_{51}+I_{52}.
\end{array}
$$
One has
$$
\begin{array}{ll}
\di I_{51}&\di \leq C\int_0^t\|(\bar u_{xx},\bar
u_x^2)\|_{L^2(\mathbf{R})}\|(\r^\a-\bar\r^\a){\bf
1}|_{\{0\leq\r\leq2\r_+\}}\|_{L^2(\mathbf{R})}dxd\tau\\
&\di \leq C\sup_{t\in[0,T]}\|(\r^\a-\bar\r^\a){\bf
1}|_{\{0\leq\r\leq2\r_+\}}\|_{L^2(\mb{R})}\int_0^t\|(\bar
u_{xx},\bar
u_x^2)\|_{L^2(\mathbf{R})}d\tau\\
&\di \leq
C\sup_{t\in[0,T]}\|\r\Psi(\r,\bar\r)\|^{\f12}_{L^1(\mathbf{R})},
\end{array}
$$
due to the fact that
$$
\begin{array}{ll}
\di \f{(\r^\a-\bar\r^\a)^2{\bf
1}|_{\{0\leq\r\leq2\r_+\}}}{\r\Psi(\r,\bar\r)}&\di
=\f{(\g-1)(\r^\a-\bar\r^\a)^2{\bf
1}|_{\{0\leq\r\leq2\r_+\}}}{\r^\g-\bar\r^\g-\g\bar\r^{\g-1}(\r-\bar\r)}\\[3mm]
&\di \leq
C\f{(\r^{2\a-2}+\bar\r^{2\a-2})(\r-\bar\r)^2}{(\r-\bar\r)^2}{\bf
1}|_{\{0\leq\r\leq2\r_+\}}\\[3mm]
&\di \leq C,
\end{array}
$$
where one has used $1<\g\leq2$ and $\a\geq1$.

On the other hand,
$$
\begin{array}{ll}
\di I_{52}&\di \leq C\int_0^t\|(\bar u_{xx},\bar
u_x^2)\|_{L^\i(\mathbf{R})}\|(\r^\a-\bar\r^\a)|{\bf
1}|_{\{\r\geq2\r_+\}}\|_{L^1(\mathbf{R})}dxd\tau\\
&\di \leq C\sup_{t\in[0,T]}\|(\r^\a-\bar\r^\a){\bf
1}|_{\{\r\geq2\r_+\}}\|_{L^1(\mb{R})}\int_0^t\|(\bar u_{xx},\bar
u_x^2)\|_{L^\i(\mathbf{R})}d\tau\\
&\di \leq C
\eta^{\f2{4q+1}}\sup_{t\in[0,T]}\|\r\Psi(\r,\bar\r)\|_{L^1(\mathbf{R})},
\end{array}
$$
due to the fact that
$$
\begin{array}{ll}
\di \lim_{\r\rightarrow+\i}\f{|\r^\a-\bar\r^\a|{\bf
1}|_{\{\r\geq2\r_+\}}}{\r\Psi(\r,\bar\r)}&\di
=\lim_{\r\rightarrow+\i}\f{(\g-1)|\r^\a-\bar\r^\a|{\bf
1}|_{\{\r\geq2\r_+\}}}{\r^\g-\bar\r^\g-\g\bar\r^{\g-1}(\r-\bar\r)}\\
&\di \leq C,\qquad {\rm if}~~\a\leq \g.
\end{array}
$$
Finally, we estimate the terms $I_4$ and $I_6$. As for $I_3$, one
has
$$
\begin{array}{ll}
\di (\bar\r^{\f{\t+\g-1}{2}})_{xx}&\di
=(\f{\t+\g-1}{2}\bar\r^{\f{\t+\g-3}{2}}\bar\r_x)_x=\f{\t+\g-1}{2\sqrt\g}(\bar\r^{\f\t2}\bar
u_x)_x\\
&\di =\f{\t+\g-1}{2\sqrt\g}\bar\r^{\f\t2}\bar
u_{xx}+\f{\t(\t+\g-1)}{4\g}\bar\r^{\f{\t+1-\g}{2}}\bar u_x^2,
\end{array}
$$
therefore,
\begin{equation}
\begin{array}{ll}
I_4&\di\leq C\v\int_0^t\int_{\mathbf{R}}|(\bar
u_{xx},\bar\r^{\f{\t+1-\g}{2}}\bar
u_x^2)||(\r^{\f{\t+1-\g}{2}}-\bar\r^{\f{\t+1-\g}{2}})|dxdt\\
&\di \leq C\v\nu^{\f{\t-1}{2}}\int_0^t\int_{\mathbf{R}}|(\bar
u_{xx},\bar
u_x^2)||(\r^{\f{\t+1-\g}{2}}-\bar\r^{\f{\t+1-\g}{2}})|dxdt\\
&\di = C\v\nu^{\f{\t-1}{2}}\int_0^t\int_{\mathbf{R}}|(\bar
u_{xx},\bar
u_x^2)||(\r^{\f{\t+1-\g}{2}}-\bar\r^{\f{\t+1-\g}{2}})|[{\bf
1}|_{\{0\leq\r<\f{\nu}2\}}+{\bf
1}|_{\{\f{\nu}2\leq\r\leq2\r_+\}}+{\bf 1}|_{\{\r>2\r_+\}}]dxdt\\
&\di :=I_{41}+I_{42}+I_{43}.
\end{array}
\end{equation}
Recall the following useful fact that for any given $C>0,$ there
exists a constant $\b\in[0,1]$ such that
\begin{equation}
\begin{array}{ll}
\di \r\Psi(\r,\bar\r){\bf 1}|_{\{0\leq\r\leq C\}}&\di =\f{1}{\g-1}\Big[\r^\g-\bar\r^\g-\g\bar\r^{\g-1}(\r-\bar\r)\Big]{\bf 1}|_{\{0\leq\r\leq C\}}\\
&\di =\big[\b\r+(1-\b)\bar\r\big]^{\g-2}(\r-\bar\r)^2{\bf 1}|_{\{0\leq\r\leq C\}}\\
 &\di \geq \max\{\r_+,C\}^{\g-2}(\r-\bar\r)^2,
\end{array}  \label{fact1}
\end{equation}
provided  that $1<\g\leq 2.$

One can compute that
\begin{equation}
\begin{array}{ll}
I_{41}&\di \leq C\v\nu^{\f{\t-1}{2}}\int_0^t\int_{\mathbf{R}}\|(\bar
u_{xx},\bar
u_x^2)\|_{L^2(\mathbf{R})}\|(\r^{\f{\t+\g-1}{2}}-\bar\r^{\f{\t+\g-1}{2}}){\bf
1}_{\{0\leq\r\leq\f\nu2\}}\|_{L^2(\mathbf{R})}dxd\tau\\[3mm]
&\di\leq
C\v\nu^{\f{\t-1}{2}}\sup_{t\in[0,T]}\|(\r^{\f{\t+\g-1}{2}}-\bar\r^{\f{\t+\g-1}{2}}){\bf
1}_{\{0\leq\r<\f\nu2\}}\|_{L^2(\mb{R})}\int_0^t\|(\bar u_{xx},\bar
u_x^2)\|_{L^2(\mathbf{R})}d\tau\\[3mm]
&\di \leq
C\v\nu^{\f32(\t-1)}\sup_{t\in[0,T]}\|\r\Psi(\r,\bar\r)\|^{\f12}_{L^1(\mathbf{R})},
\end{array}\label{I41}
\end{equation}
where one has used the facts that
$$
\f{(\r^{\f{\t+\g-1}{2}}-\bar\r^{\f{\t+\g-1}{2}})^2{\bf
1}_{\{0\leq\r<\f\nu2\}}}{\r\Psi(\r,\bar\r)}\leq
\f{(\r^{\f{\t+\g-1}{2}}-\bar\r^{\f{\t+\g-1}{2}})^2{\bf
1}_{\{0\leq\r<\f\nu2\}}}{(\r_+)^{\g-2}(\r-\bar\r)^2}, ~~{\rm
if}~~\nu\ll1,
$$
and  the function
$$
\f{(\r^{\f{\t+\g-1}{2}}-\bar\r^{\f{\t+\g-1}{2}})^2}{(\r-\bar\r)^2}
$$
is monotone decreasing in $\r\in[0,\f\nu2],$ that is,
$$
\f{(\r^{\f{\t+\g-1}{2}}-\bar\r^{\f{\t+\g-1}{2}})^2{\bf
1}_{\{0\leq\r<\f\nu2\}}}{(\r-\bar\r)^2}\leq
\lim_{\r\rightarrow0}\f{(\r^{\f{\t+\g-1}{2}}-\bar\r^{\f{\t+\g-1}{2}})^2}{(\r-\bar\r)^2}=\bar\r^{\t+\g-3}\leq
\nu^{\t-1}.
$$
Moreover, it holds that
\begin{equation}
\begin{array}{ll}
I_{42}&\di \leq C\v\nu^{\f{\t-1}{2}}\int_0^t\int_{\mathbf{R}}\|(\bar
u_{xx},\bar
u_x^2)\|_{L^2(\mathbf{R})}\|(\r^{\f{\t+\g-1}{2}}-\bar\r^{\f{\t+\g-1}{2}}){\bf
1}_{\{\f\nu2\leq\r\leq2\r_+\}}\|_{L^2(\mathbf{R})}dxd\tau\\[3mm]
&\di\leq
C\v\nu^{\f{\t-1}{2}}\sup_{t\in[0,T]}\|(\r^{\f{\t+\g-1}{2}}-\bar\r^{\f{\t+\g-1}{2}}){\bf
1}_{\{\f\nu2\leq\r\leq2\r_+\}}\|_{L^2(\mb{R})}\int_0^t\|(\bar
u_{xx},\bar
u_x^2)\|_{L^2(\mathbf{R})}d\tau\\[3mm]
&\di \leq
C\v\nu^{\f32(\t-1)}\sup_{t\in[0,T]}\|\r-\bar\r\|_{L^2(\mathbf{R})}\\
&\di \stackrel{\eqref{fact1}}{\leq}
C\v\nu^{\f32(\t-1)}\sup_{t\in[0,T]}\|\r\Psi(\r,\bar\r)\|^{\f12}_{L^1(\mathbf{R})},
\end{array}\label{I42}
\end{equation}
where in the third inequality one has used the fact that there
exists a constant $\b\in[0,1],$
$$
\begin{array}{ll}
\di \f{(\r^{\f{\t+\g-1}{2}}-\bar\r^{\f{\t+\g-1}{2}})^2{\bf
1}_{\{\f\nu2\leq\r\leq2\r_+\}}}{(\r-\bar\r)^2}&\di=(\t+\g-1)^2
\big[\b\r+(1-\b)\bar\r\big]^{\t+\g-3}{\bf
1}_{\{\f\nu2\leq\r\leq2\r_+\}}\\
&\di \leq (\t+\g-1)^2(\f\nu2)^{\t-1}.
\end{array}
$$
And then
\begin{equation}
\begin{array}{ll}
I_{43}&\di \leq
C\v\nu^{\f{\t-1}{2}}\int_0^t\int_{\mathbf{R}}\|(\r^{\f{\t+\g-1}{2}}-\bar\r^{\f{\t+\g-1}{2}}){\bf
1}_{\{\r>2\r_+\}}\|_{L^1(\mathbf{R})}\|(\bar u_{xx},\bar
u_x^2)\|_{L^\i(\mathbf{R})}dxd\tau\\[3mm]
&\di\leq
C\v\nu^{\f{\t-1}{2}}\sup_{t\in[0,T]}\|(\r^{\f{\t+\g-1}{2}}-\bar\r^{\f{\t+\g-1}{2}}){\bf
1}_{\{\r>2\r_+\}}\|_{L^1(\mb{R})}\int_0^t\|(\bar u_{xx},\bar
u_x^2)\|_{L^\i(\mathbf{R})}d\tau\\[3mm]
&\di \leq
C\v\nu^{\f{\t-1}2}\sup_{t\in[0,T]}\|\r\Psi(\r,\bar\r)\|_{L^1(\mathbf{R})},
\end{array}\label{I43}
\end{equation}
since
$$
\lim_{\r\ra+\i}\f{(\r^{\f{\t+\g-1}{2}}-\bar\r^{\f{\t+\g-1}{2}}){\bf
1}_{\{\r>2\r_+\}}}{\r\Psi(\r,\bar\r)}=0.
$$
In summary, by combining \eqref{I41}, \eqref{I42} and \eqref{I43},
one can arrive at
\begin{equation}
I_4\leq
C\v\nu^{\f32(\t-1)}\Big[\sup_{t\in[0,T]}\|\r\Psi(\r,\bar\r)\|_{L^1(\mathbf{R})}+1\Big].
\end{equation}
Finally, $I_6$ can be estimated similarly as for $I_4$ and  the
details will be omitted for brevity.

 For definiteness, we take $\v\nu^{\f32(\t-1)}=\v^{\f12}$,
i.e., $\nu=\v^{\f23}$ since $\t=\f12.$ Consequently, choosing $\v$
such that $\v \ln(1+T)\leq 1$ and $\v$ suitably small and combining
all the above estimates shows that for $\a$ and $\g$ satisfying
\eqref{alpha},
\begin{equation}
\begin{array}{ll}
\di \sup_{t\in[0,T]}\int_{\mathbf{R}}\Big\{\r(u-\bar
u)^2+\Big[\Big(\f{\r^{\a-\f12}}{\a-\f12}\Big)_x\Big]^2+\v^2\big[\big(\ln\r\big)_x\big]^2+\r\Psi(\r,\bar\r)\Big\}(x,t)dx\\[3mm]
\di +\int_0^T\int_{\mathbf{R}}\Big\{\bar
u_x\Big[p(\r)-p(\bar\r)-p^\prime(\bar\r)(\r-\bar\r)\Big]+\r(u-\bar
u)^2\bar u_x+(\r^\a+\v\r^\t)\Big[(u-\bar u)_x\Big]^2\\
\di \qquad\qquad~~
+\Big[(\r^{\f{\a+\g-1}{2}}-\bar\r^{\f{\a+\g-1}{2}})_x\Big]^2
+\v\Big[(\r^{\f{\t+\g-1}{2}}-\bar\r^{\f{\t+\g-1}{2}})_x\Big]^2\Big\}(x,t)dxdt
\leq C.
\end{array}
\end{equation}
Thus Lemma \ref{lemma1} is proved.  $\hfill\Box$

The following lemma is the key point to get the existence of the
approximate solution $(\r_{\v,\nu},u_{\v,\nu})(t,x)$ with
$\nu=\v^{\f23}$.

\begin{Lemma}\label{lemma2} There exist an absolutely constant $C$ and a
positive constant $C(\v,\nu,T)$ depending on $\v,\nu$ and $T$  such
that
\begin{eqnarray}
0<C(\v,\nu,T)\le\rho_{\v,\nu}\le C.\label{D1}
\end{eqnarray}
\end{Lemma}

{\em Proof:}
 From the Gagliardo-Nirenberg inequality:
$$
\|f\|_{L^\i(\mathbf{R})}\leq
C\|f_x\|_{L^2(\mathbf{R})}^\b\|f\|_{L^p(\mathbf{R})}^{1-\b},
$$
where $0<\b<1$, $1\leq p<\i$ to be determined, and $\b$, $p$ satisfy
$$
\f{\b}{2}=\f{1-\b}{p},
$$
we have
\begin{equation}
\begin{array}{ll}
\di \sup_{t\in[0,T]}\|\r^{\a-\f12}-\bar\r^{\a-\f12}\|_{L^\i(\mathbf{R})}\\
\di\leq
C\sup_{t\in[0,T]}\|(\r^{\a-\f12}-\bar\r^{\a-\f12})_x\|_{L^2(\mathbf{R})}^\b\sup_{t\in[0,T]}\|\r^{\a-\f12}-\bar\r^{\a-\f12}\|_{L^p(\mathbf{R})}^{1-\b}\\
\di\leq
C\sup_{t\in[0,T]}\Big[\|(\r^{\a-\f12})_x\|_{L^2(\mathbf{R})}^\b+\|(\bar\r^{\a-\f12})_x\|_{L^2(\mathbf{R})}^\b\Big]
\sup_{t\in[0,T]}\|\r^{\a-\f12}-\bar\r^{\a-\f12}\|_{L^p(\mathbf{R})}^{1-\b}
\\\di\leq
C\sup_{t\in[0,T]}\|\r^{\a-\f12}-\bar\r^{\a-\f12}\|_{L^p(\mathbf{R})}^{1-\b},
\end{array}\label{upper}
\end{equation}
due to the fact that
$$
\|(\bar\r^{\a-\f12})_x\|_{L^2(\mathbf{R})}=\|(\a-\f12)\bar\r^{\a-\f\g2}\bar
u_x\|_{L^2(\mathbf{R})}\leq C\|\bar u_x\|_{L^2(\mathbf{R})}\leq C.
$$
Since
\begin{equation}
\lim_{\r\rightarrow0}\f{|\r^{\a-\f12}-\bar\r^{\a-\f12}|^{\f{2\g}{2\a-1}}}{\r\Psi(\r,\bar\r)}=1,
\end{equation}
there exists a positive constant $\d^1_{\f12}$, such that if $0\leq
\r\leq \d^1_{\f12}$, then
\begin{equation}
|\f{|\r^{\a-\f12}-\bar\r^{\a-\f12}|^{\f{2\g}{2\a-1}}}{\r\Psi(\r,\bar\r)}-1|<\f12,
\end{equation}
thus, for any $\bar\r\geq0,$
\begin{equation}
\f12\r\Psi(\r,\bar\r)\leq
|\r^{\a-\f12}-\bar\r^{\a-\f12}|^{\f{2\g}{2\a-1}}\leq
\f32\r\Psi(\r,\bar\r),~~{\rm if}~~0\leq \r\leq
\d^1_{\f12}.\label{p0}
\end{equation}
Similarly, there exists a positive constant $\bar\d^1_{\f12}$, such
that for any $\r\geq 0,$
\begin{equation}
\f12\r\Psi(\r,\bar\r)\leq
|\r^{\a-\f12}-\bar\r^{\a-\f12}|^{\f{2\g}{2\a-1}}\leq
\f32\r\Psi(\r,\bar\r),~~{\rm if}~~0\leq \bar\r\leq
\bar\d^1_{\f12}.\label{p1}
\end{equation}
Set
$$
p=\f{2\g}{2\a-1}\in[2,2\g].
$$
Then for such $p=\f{2\g}{2\a-1}$ and $\g\in(1,2]$, one has
\begin{equation}\label{p21}
\f{|\r^{\a-\f12}-\bar\r^{\a-\f12}|^p}{\r\Psi(\r,\bar\r)}{\bf
1}|_{\{\d^1_{\f12}\leq\r\leq2\r_+,~\bar\d^1_{\f12}\leq\bar\r\leq\r_+\}}\leq
C\f{|\r^{\a-\f12}-\bar\r^{\a-\f12}|^2}{(\r-\bar\r)^2}{\bf
1}|_{\{\d^1_{\f12}\leq\r\leq2\r_+,~\bar\d^1_{\f12}\leq\bar\r\leq\r_+\}}\leq
C,
\end{equation}
and
\begin{equation}\label{p2+}
\lim_{\r\rightarrow+\i}\f{|\r^{\a-\f12}-\bar\r^{\a-\f12}|^p}{\r\Psi(\r,\bar\r)}
=\lim_{\r\rightarrow+\i}
\f{(\g-1)|\r^{\a-\f12}-\bar\r^{\a-\f12}|^{\f{2\g}{2\a-1}}}{\r^\g-\bar\r^\g-\g\bar\r^{\g-1}(\r-\bar\r)}
=\g-1.
\end{equation}
Thus it follows from \eqref{p2+} that
\begin{equation}
|\r^{\a-\f12}-\bar\r^{\a-\f12}|^p\leq C\r\Psi(\r,\bar\r),~~{\rm
if}~\r\geq2\r_+.\label{p2}
\end{equation}
Collecting \eqref{p0}, \eqref{p1}, \eqref{p21} and \eqref{p2} gives
that
\begin{equation}
\begin{array}{ll}
\di\|\r^{\a-\f12}-\bar\r^{\a-\f12}\|^p_{L^p(\mathbf{R})}
=\int_{\mathbf{R}}|\r^{\a-\f12}-\bar\r^{\a-\f12}|^pdx\\[2mm]
\di =\int_{\mathbf{R}}|\r^{\a-\f12}-\bar\r^{\a-\f12}|^p\big({\bf
1}|_{\{0\leq\r\leq\d^1_{\f12}\}}+{\bf
1}|_{\{\d^1_{\f12}\leq\r\leq2\r_+,~0\leq\bar\r\leq\bar\d^1_{\f12}\}}\\
\qquad\qquad+{\bf
1}|_{\{\d^1_{\f12}\leq\r\leq2\r_+,~\bar\d^1_{\f12}\leq\bar\r\leq\r_+\}}+{\bf
1}|_{\{\r\geq2\r_+\}}\big)dx\\
\di \leq C\int_{\mathbf{R}}\r\Psi(\r,\bar\r)dx \leq C.
\end{array}\label{p3}
\end{equation}
Substituting \eqref{p3} into \eqref{upper} yields the uniform upper
bound for $\r_{\v,\nu}(t,x)$.

Next we derive a lower bound for $\r_{\v,\nu}(t,x)$. Since
$\lim_{\r\rightarrow0}\r\Psi(\r,\bar\r)=\bar\r^{\g}\geq\nu^\g$, then
$\r\Psi(\r,\bar\r)$ is bounded away from 0 on $[0,\f12\bar\r]$. Thus
one can deduce from the bound on $\r\Psi(\r,\bar\r)$ in
$L^\i(0,T;L^1(\mathbf{R}))$ that there exists a constant
$C_1=C_1(\nu,T)>0$, such that for all $t\in[0,T]$,
$$
{\rm meas} \{x\in\mathbf{R}|\r(x,t)\leq \f12\bar\r(x,t)\}\leq
\f{1}{\inf_{\r\in[0,\f12\bar\r]}\r\Psi(\r,\bar\r)}\int_{\{x\in\mathbf{R}|\r(x,t)\leq
\f12\bar\r(x,t)\}}\r\Psi(\r,\bar\r)(x,t)dx\leq C_1.
$$
Therefore, for every $x_0\in\mathbf{R}$, there exists $M=M(\nu,T)>0$
large enough, such that
$$
\begin{array}{ll}
\di \int_{|x-x_0|\leq M}\r_{\v,\nu}(x,t)dx&\di
\geq\int_{\{|x-x_0|\leq
M\}\cap\{x\in\mathbf{R}|\r_{\v,\nu}(x,t)> \f12\bar\r(x,t)\}}\r_{\v,\nu}(x,t)dx\\
&\di \geq \f12\inf_{(x,t)}\bar\r(x,t){\rm meas}\Big\{\{|x-x_0|\leq
M\}\cap\{x\in\mathbf{R}|\r_{\v,\nu}(x,t)> \f12\bar\r(x,t)\}\Big\}\\
&\di =\f\nu2{\rm meas}\Big\{\{|x-x_0|\leq
M\}\cap\{x\in\mathbf{R}|\r_{\v,\nu}(x,t)\leq \f12\bar\r(x,t)\}^{c}\Big\}\\
&\di \geq \f\nu2(2M-C_1)>0,
\end{array}
$$
for all $t\in [0,T]$.

Due to the continuity of $\r_{\v,\nu}$, there exists $x_1\in
[x_0-M,x_0+M]$ such that
$$
\r_{\v,\nu}(x_1,t)=\f{1}{2M}\int_{|x-x_0|\leq
M}\r_{\v,\nu}(x,t)dx\geq \f\nu{4M} (2M-C_1).
$$
Thus,
$$
\begin{array}{ll}
\di |\ln\r_{\v,\nu}(x_0,t)|&\di =|\ln\r_{\v,\nu}(x_1,t)+\int_{x_1}^{x_0}(\ln\r_{\v,\nu})_x(x,t)dx|\\
&\di \leq |\ln\r_{\v,\nu}(x_1,t)|+\|(\ln\r_{\v,\nu})_x(\cdot,t)\|_{L^2(\mathbf{R})}|x_1-x_0|^{\f12}\\
&\di \leq C(\v,\nu,M,T)+C_\v M^{\f12}.
\end{array}
$$
Consequently, we can get that there exists a positive constant
$C(\v,\nu,T)$ such that
$$
\r_{\v,\nu}(x_0,t)\geq C(\v,\nu,T),
$$
for any $x_0\in\mathbf{R}$ and $t\in[0,T]$. $\hfill\Box$

With the lower and upper bounds on $\r_{\v,\nu}$, we can get the
existence of the approximate solution
$(\r_{\v,\nu},u_{\v,\nu})(t,x)$ by a similar argument as in
\cite{MV2}. In order to pass the limit $\v\ra0$ with
$\nu=\v^{\f23}$, we need the following higher estimates on the
momentum.

\begin{Lemma}\label{lemma3} There exists a positive constant $C(T)$
independent of $\v,\nu$, such that
$$
\sup_{t\in[0,T]}\int_{\mathbf{R}}\r_{\v,\nu}|u_{\v,\nu}-\bar
u_\nu|^{3}(t,x)dx+\int_0^T\int_{\mathbf{R}}\r_{\v,\nu}[(u_{\v,\nu}-\bar
u_{\nu})_x]^2|u_{\v,\nu}-\bar u| dxdt\leq C(T).
$$
\end{Lemma}
The proof of Lemma \ref{lemma3} can be done along the same line as
in our previous paper \cite{JWX} and we omit the details for
brevity.

Now with these uniform in $\v,\nu$ estimates at hand, we can pass
the limit process $\v\rightarrow0$ with $\nu=\v^{\f23}$, obtain the
existence of the weak solution $(\r,u)(t,x)$, and get the uniform in
time estimates in Theorem \ref{w-e}.

\section{Asymptotic behavior of weak solutions}

In this section, we will study the asymptotic behavior of a given
weak solution $(\r,u)(t,x)$ to the Cauchy problem
\eqref{(1.1)}-\eqref{(1.3)} in the sense of Definition \ref{def}
satisfying  \eqref{P12-1} and \eqref{P12-3}.  We assume that the
solution is smooth enough. The rigorous proof can be obtained by
using the usual regularization procedure.

{\bf \underline{Proof of Theorem \ref{asymptotic-b}}:} For any $s\ge
1$, by the uniform upper bound of $\r$, it holds that
\begin{eqnarray*}
|\rho^s-\bar\rho^s|^2\le C|\rho-\bar\rho|^2.
\end{eqnarray*}
Hence it follows from \eqref{P12-3} and \eqref{fact1} that
\begin{eqnarray}
\int_{\mathbf{R}} |\rho^s-\bar\rho^s|^2 dx \le C\int_{\mathbf{R}}
(\rho-\bar\rho)^2 dx\le C.\label{Jan6-2}
\end{eqnarray}
Similarly,
\begin{eqnarray}
\int_{\mathbf{R}} |\rho^s-\bar\rho^s|^{2\l} dx \le
C\int_{\mathbf{R}} |\rho-\bar\rho|^{2\l} dx\leq C\int_{\mathbf{R}}
(\rho-\bar\rho)^2 dx\le C,\label{Jan6-4}
\end{eqnarray}
for any $\lambda\ge 1$.

Set $b=\f{\a+\g-1}{2}$. Then one gets from \eqref{P12-3} and
\eqref{fact1} that
$$
\int_0^t\int_{\mathbf{R}}\big\{[(\r^b-\bar\r^b)_x]^2+\bar
u_x(\r-\bar\r)^2\big\}dxd\tau\leq C.
$$
For $s>b+1$, it holds that
\begin{eqnarray*}
&&(\rho^s-\bar\rho^s)^2(t,x)=\int_{-\infty}^x
[(\rho^s-\bar\rho^s)^2]_x
dx\\
&&=2\int_{-\infty}^x (\rho^s-\bar\rho^s)(\rho^s-\bar\r^s)_x dx\\
&&=2s\int_{-\infty}^x
(\rho^s-\bar\rho^s)[(\rho^b-\bar\r^b)_x\rho^{s-1}+(\bar\r^b)_x(\r^{s-b}-\bar\r^{s-b})] dx\\
&&\le C
\|\rho^s-\bar\rho^s\|_{L^2(\mathbf{R})}\|(\rho^b-\bar\r^b)_x\|_{L^2(\mathbf{R})}+C\int_\mathbf{R}\bar\r^{\f\a2}\bar u_x(\r^s-\bar\r^s)(\r^{s-b}-\bar\r^{s-b})dx\\
&&\leq
\|\rho^s-\bar\rho^s\|_{L^2(\mathbf{R})}\|(\rho^b-\bar\r^b)_x\|_{L^2(\mathbf{R})}+C\int_\mathbf{R}\bar
u_x(\r-\bar\r)^2dx.
\end{eqnarray*}
Consequently,
\begin{eqnarray}
&&\int_0^t\sup_{x\in \mathbf{R}}(\rho^s-\bar\rho^s)^4
dt\nonumber\\&&\le
C\sup_{t\in[0,T]}\|\rho^s-\bar\rho^s\|^2_{L^2(\mathbf{R})}\int_0^t
\|(\rho^b-\bar\r^b)_x\|_{L^2(\mathbf{R})}^2 dt\\
&&\quad+C\sup_{t\in[0,T]}\|\r-\bar\r\|^2_{L^2(\mathbf{R})}\int_0^t\int_\mathbf{R}\bar
u_x(\r-\bar\r)^2dxd\tau\nonumber\\
&&\le C\nonumber.
\end{eqnarray}
Moreover, applying \eqref{Jan6-4} leads to
\begin{equation}\begin{array}{l} \di \int_0^t\int_\mathbf{R}
(\rho^s-\bar\rho^s)^4(\rho^s-\bar\rho^s)^{2l}
dxdt\\
\di\leq \int_0^t\Big[\sup_{x\in\mathbf{R}}(\rho^s-\bar\rho^s)^4\int_{\mathbf{R}}(\rho^s-\bar\rho^s)^{2l}dx\Big]dt\\
 \di \le \sup_{t\in[0,T]} \int_\mathbf{R} (\rho^s-\bar\rho^s)^{2l}
dx\int_0^t \sup_{x\in \mathbf{R}}(\rho^s-\bar\rho^s)^4 dt\le
C,~~~\forall l\ge 1.
\end{array}\label{(4.12)}\end{equation}
Set
$$ f(t)=\int_{\mathbf{R}}
(\rho^s-\bar\rho^s)^{4+2l} dx.
 $$
Then $$ f(t)\in
 L^1(0,\infty)\cap L^\infty(0,\infty)$$ due to \eqref{Jan6-4} and \eqref{(4.12)}.

  Furthermore, direct
 calculations show that
\begin{equation}
\begin{array}{ll}
\di \frac{d}{dt} f(t) =(4+2l)s\int_\mathbf{R}
(\rho^s-\bar\rho^s)^{3+2l}
(\rho^{s-1}\rho_t-\bar\r^{s-1}\bar\r_t) dx\\
~~\di =-(4+2l)s\int_\mathbf{R}
(\rho^s-\bar\rho^s)^{3+2l}[\rho^{s-1}(\rho
u)_x-\bar\r^{s-1}(\bar\r\bar u)_x]
dx\\
~~\di =(4+2l)(3+2l)s\int_\mathbf{R}
(\rho^s-\bar\rho^s)^{2+2l}(\rho^s-\bar\r^s)_x(\rho^{s-1}\rho
u-\bar\r^{s-1}\bar\r\bar u)
dx\\
~~\di \quad+(4+2l)s(s-1)\int_\mathbf{R}
(\rho^s-\bar\rho^s)^{3+2l}(\rho^{s-2}\rho_x\rho
u-\bar\rho^{s-2}\bar\rho_x\bar\r \bar u)
dx\\
\di =(4+2l)(3+2l)s\int_\mathbf{R}
(\rho^s-\bar\rho^s)^{2+2l}(\rho^s-\bar\r^s)_x\rho^{s}(u-\bar u)
dx\\
\quad\di +(4+2l)(3+2l)s\int_\mathbf{R}
(\rho^s-\bar\rho^s)^{2+2l}(\rho^s-\bar\r^s)_x\bar
u(\rho^{s}-\bar\r^{s})
dx\\
\di \quad+(4+2l)s(s-1)\int_\mathbf{R}
(\rho^s-\bar\rho^s)^{3+2l}\rho^{s-1}\rho_x
(u-\bar u) dx\\
\di \quad+(4+2l)s(s-1)\int_\mathbf{R} (\rho^s-\bar\rho^s)^{3+2l}\bar
u(\rho^{s-1}\rho_x
-\bar\rho^{s-1}\bar\rho_x) dx\\
\di :=J_1(t)+J_2(t)+J_3(t)+J_4(t).
\end{array}
\end{equation}
Now we claim that $J_i(t)\in L^2(0,+\i),~(i=1,2,3,4)$. In fact,
\begin{equation}
\begin{array}{ll}
\di J_1(t)&\di=\f{(4+2l)(3+2l)s^2}{b}\int_\mathbf{R}
(\rho^s-\bar\rho^s)^{2+2l}\rho^{s}(u-\bar
u)\Big[\r^{s-b}(\rho^b-\bar\r^b)_x+(\r^{s-b}-\bar\r^{s-b})(\bar\r^b)_x\Big]
dx\\
&\di =\f{(4+2l)(3+2l)s^2}{b}\int_\mathbf{R}
(\rho^s-\bar\rho^s)^{2+2l}\sqrt{\r}(u-\bar
u)\r^{2s-b-\f12}(\rho^b-\bar\r^b)_xdx\\
&\di\quad +\f{(4+2l)(3+2l)s^2}{b}\int_\mathbf{R}
(\rho^s-\bar\rho^s)^{2+2l}\sqrt{\r}(u-\bar
u)\r^{s-\f12}(\r^{s-b}-\bar\r^{s-b})(\bar\r^b)_x dx\\
&\leq C\|\sqrt{\r}(u-\bar
u)\|_{L^2(\mathbf{R})}\|(\rho^b-\bar\r^b)_x\|_{L^2(\mathbf{R})}\\
&\di\quad +C\|\sqrt{\r}(u-\bar
u)\|_{L^2(\mathbf{R})}\|(\rho^s-\bar\rho^s)^{2+2l}(\r^{s-b}-\bar\r^{s-b})(\bar\r^b)_x
\|_{L^2(\mathbf{R})}\\
&\leq C\|\sqrt{\r}(u-\bar
u)\|_{L^2(\mathbf{R})}\|(\rho^b-\bar\r^b)_x\|_{L^2(\mathbf{R})}+C\|\sqrt{\r}(u-\bar
u)\|_{L^2(\mathbf{R})}\|\bar u_x(\r-\bar\r) \|_{L^2(\mathbf{R})}.
\end{array}
\end{equation}
Thus,
\begin{equation}
\begin{array}{ll}
\di \int_0^t |J_1(t)|^2dt&\di \leq
C\sup_{t\in[0,T]}\|\sqrt{\r}(u-\bar
u)\|^2_{L^2(\mathbf{R})}\int_0^t\|(\rho^b-\bar\r^b)_x\|^2_{L^2(\mathbf{R})}dt\\
&\di\quad +C\sup_{t\in[0,T]}\|\sqrt{\r}(u-\bar
u)\|^2_{L^2(\mathbf{R})}\int_0^t\|\bar u_x(\r-\bar\r)
\|_{L^2(\mathbf{R})}^2dt\\
&\di \leq C\sup_{t\in[0,T]}\|\sqrt{\r}(u-\bar
u)\|^2_{L^2(\mathbf{R})}\int_0^t\|(\rho^b-\bar\r^b)_x\|^2_{L^2(\mathbf{R})}dt\\
&\di\quad +C\sup_{t\in[0,T]}\|\sqrt{\r}(u-\bar
u)\|^2_{L^2(\mathbf{R})}\int_0^t\int_{\mathbf{R}}\bar
u_x(\r-\bar\r)^2 dxdt\\
&\di \leq C.
\end{array}
\end{equation}
The fact that $J_i(t)\in L^2(0,+\i),~(i=2,3,4)$ can be shown
similarly.

 Thus
$$
\f{d}{dt}f(t)\in L^2(0,+\i).
$$
Combining the fact that $f(t)\in
 L^1(0,+\infty)\cap L^\infty(0,+\infty)$, one has
\begin{eqnarray}
f(t)\to 0, \ \ t\to+\infty.
\end{eqnarray}
Let $m\ge 1$ be any real number to be determined later. One has
\begin{equation}
\begin{array}{ll}
\di |(\rho^s-\bar\rho^s)^m| =|\int_{-\infty}^x
[(\rho^s-\bar\rho^s)^m]_x dx|\\
\quad\di =|m\int_{-\infty}^x
(\rho^s-\bar\rho^s)^{m-1}(\rho^s-\bar\r^s)_x dx|\\
\quad\di =|ms\int_{-\infty}^x
(\rho^s-\bar\rho^s)^{m-1}\Big[\f{1}{\a-\f12}(\rho^{\a-\f12})_x\r^{s-\a+\f12}-\bar\r^{s-1}\bar\r_x \Big]dx|\\
\quad \di \leq
C\|(\rho^s-\bar\rho^s)^{m-1}\|_{L^2(\mathbf{R})}\Big[\|(\rho^{\a-\f12})_x\|_{L^2(\mathbf{R})}+\|\bar\r^{\f{3-\g}{2}}\bar u_x\|_{L^2(\mathbf{R})}\Big]\\
\quad\di \le C\|(\rho^s-\bar\rho^s)^{m-1}\|_{L^2(\mathbf{R})}.
\end{array}
\end{equation}
Choosing $2(m-1)=4+2l$ yields
\begin{eqnarray}
\sup_{x\in \mathbf{R}}|(\rho^s-\bar\rho^s)^m|\le C f^\frac12(t)\to
0,~~\mbox{as}~~ t\to+\infty.
\end{eqnarray}
Therefore,
$$
\lim_{t\to+\infty}\sup_{x\in\mathbf{R}}|\rho^s-\bar\rho^s|=0.
$$
Now we prove that
$\di\lim_{t\to+\infty}\sup_{x\in\mathbf{R}}|\rho-\bar\rho|=0$.

Since
$$
\lim_{\bar\r\rightarrow0+}\f{|\r-\bar\r|^s}{|\r^s-\bar\r^s|}=1,
\quad{\rm uniformly ~in}~~\r,
$$
then for any $\sigma>0$, there exists $\d_\sigma>0$, such that if
$0\leq\bar\r<\d_\sigma$, then
$$
|\f{|\r-\bar\r|^s}{|\r^s-\bar\r^s|}-1|\leq\sigma.
$$
Thus, fix $\sigma=\f12$, then if $0\leq\bar\r<\d:=\d_{\f12}$, one
has for any $\r\geq0,$
\begin{equation}
|\r-\bar\r|^s\leq \f32|\r^s-\bar\r^s|. \label{relation}
\end{equation}
Now
$$
\begin{array}{ll}
\di |\rho-\bar\rho|^s&\di =|\rho-\bar\rho|^s({\bf
1}_{\{0\leq\bar\rho<\d\}}+{\bf
1}_{\{\bar\r\geq\d,0\leq\rho<\f{\bar\r}{2}\}}+{\bf
1}_{\{\bar\r\geq\d,\r>\f{\bar\r}{2}\}})\\
&\di\leq \f32|\rho^s-\bar\rho^s|{\bf
1}_{\{0\leq\bar\rho<\d\}}+C_\d|\rho^s-\bar\rho^s|{\bf
1}_{\{\bar\r\geq\d,0\leq\rho<\f{\bar\r}{2}\}}+C_\d|\rho^s-\bar\rho^s|^s{\bf
1}_{\{\bar\r\geq\d,\rho>\f{\bar\r}{2}\}}.
\end{array}
$$
Therefore,
\begin{eqnarray*}
&& \sup_{x\in \mathbf{R}}|\rho-\bar\rho|^s \leq
C_\d\sup_{x\in\mathbf{R}}|\rho^s-\bar\rho^s|+C_\d\sup_{x\in
\mathbf{R}}|\rho^s-\bar\rho^s|^s\to 0,
\end{eqnarray*}
 as $t\to+\infty$, which implies that
$$
 \lim_{t\to\infty}\sup_{x\in\mathbf{R}}|\rho-\bar\rho|=0.
$$
The proof of Theorem \ref{asymptotic-b} is finished.

\section{Regularity of the solution away from the vacuum}
\setcounter{equation}{0}
 In this section, we will prove Theorem \ref{regularity}, that is, we will show that away
from the vacuum region of the 2-rarefaction wave $(\r^r,u^r)(\x)$,
any weak solution $(\r,u)(t,x)$ to the Cauchy problem
\eqref{(1.1)}-\eqref{(1.3)} satisfying \eqref{P12-1} and
\eqref{P12-3} becomes regular  as stated in Theorem
\ref{regularity}.

Due to the definition of the 2-rarefaction wave $\r^r(t,x)$ in
\eqref{R2}, for any  fixed $\sigma>0$, there exist a unique
$u_\sigma$ such that $(\r,u)=(\sigma,u_\sigma)$ lies on the
2-rarefaction wave curve. In fact,
\begin{equation}
u_\sigma=\f{2\sqrt{\g}}{\g-1}\sigma^{\f{\g-1}{2}}+\Sigma_2(\r_+,u_+),\label{sigma}
\end{equation}
since 2-Riemann invariant
$\Sigma_2(\r,u)=u-\f{2\sqrt{\g}}{\g-1}\r^{\f{\g-1}{2}}$ is constant
along the 2-rarefaction wave curve. Due to the expanding property of
the 2-rarefaction wave, we have the lower bound of the density
function $\r^r(t,x)$ of 2-rarefaction wave on the right of the
straight line $x=\l^\sigma_{2} t$ with
$\l^\sigma_{2}=\l_2(\sigma,u_\sigma)=\l_2
(\r,u)|_{(\r,u)=(\sigma,u_\sigma)}$, that is,
\begin{equation}
\r^r(t,x)\geq \sigma,\quad{\rm if}\quad x\geq \l^\sigma_{2}t.
\end{equation}

Thus it follows from the asymptotic behavior \eqref{th14-1} of
$\r(t,x)$ that for $\f{\sigma}{2}>0$, there exists a large time
$T_{\sigma}$ such that if $t>T_{\sigma}$, then
\begin{equation}
\|\r(t,\cdot)-\r^r(t,\cdot)\|_{L^\i(\mathbf{R})}\leq\f{\sigma}{2}.
\end{equation}
Therefore, in the domain
\begin{equation}
\Omega_{\sigma}:=\{(t,x)|t>T_{\sigma},
x>\l^\sigma_{2}t\},\label{domain}
\end{equation}
it holds that
$$
\r(t,x)\geq\f{\sigma}{2}, \quad {\rm if}\quad (t,x)\in
\Omega_{\sigma}.
$$
So in the domain $\Omega_{\sigma}$, any vacuum states vanish and
thus
 the higher regularity of the weak solution $(\r,u)(t,x)$ can be expected as stated
in Theorem \ref{regularity}. In the following, we give the proof of
Theorem \ref{regularity}. First we establish  the local uniform
boundedness of the velocity  $u(x,t)$ by the De Giorgi-Moser
iteration method. To this end, we rewrite the momentum equation as
\begin{equation}
u_t+uu_x+\g\r^{\g-2}\r_x=\r^{\a-1}u_{xx}+\a\r^{\a-2}\r_xu_x,\label{vel}
\end{equation}

For any $(t_*,x_*)\in \Omega_{\sigma}$, and for any $r,s>0$ such
that $Q_{r,s}^*:=B_r(x_*)\times(t_*,t_*+s]\subset\Omega_{\sigma}$,
and for any test function $\z(x,t)\in
\mathring{W}_2^{1,1}(Q_{r,s}^*)$ satisfying $0\leq\z\leq 1$, one can
get from the uniform estimates in \eqref{P12-1} and \eqref{P12-3}
that
\begin{equation}
\sup_{t_*\leq t\leq t_*+s}\int_{B_r(x_*)}\Big[(u-\bar
u)^2+\r_x^2+(\r-\bar\r)^2\Big](x,t)dx+\int_{t_*}^t\int_{B_r(x_*)}\Big[(u-\bar
u)_x^2+(\r-\bar\r)_x^2\Big]dxdt\leq C.
 \label{UE1}
\end{equation}
It follows from the construction of the rarefaction wave
$(\bar\r,\bar u)(x,t)$ that
\begin{equation}
\sup_{t_*\leq t\leq
t_*+s}\int_{B_r(x_*)}(u^2+\r_x^2+\r^2)(x,t)dx+\int_{t_*}^t\int_{B_r(x_*)}(u
_x^2+\r_x^2)dxdt\leq C.
 \label{UE}
\end{equation}
Multiplying the equation \eqref{vel} by $\z^2(u-k)_+$ for any
$k\in\mathbf{R}$ and integrating the resulted equation over
$B_r(x_*)\times(t_*,t]$ for $ t\in(t_*,t_*+s]$, one arives
\begin{equation}
\begin{array}{ll}
\di \f12\int_{B_r(x_*)}\z^2(u-k)_+^2(x,t)dx+\int_{t_*}^t\int_{B_r(x_*)}\r^{\a-1}[\z(u-k)_+]_x^2dxdt\\
\di
=\f12\int_{B_r(x_*)}\z^2(u-k)_+^2(x,t_*)dx+\int_{t_*}^t\int_{B_r(x_*)}\Big\{\z\z_t(u-k)_+^2+\f23\z\z_x[(u-k)_+]^3\\\di
\qquad~-\g\r^{\g-2}\r_x\z^2(u-k)_+
+\r^{\a-1}\z_x^2(u-k)_+^2+\r^{\a-2}\r_xu_x\z^2(u-k)_+\Big\}dxdt\\
\di \leq\f12\int_{B_r(x_*)}\z^2(u-k)_+^2(x,t_*)dx+\f18\int_{t_*}^t\int_{B_r(x_*)}\r^{\a-1}[\z(u-k)_+]_x^2dxdt\\
\di\quad +
C\int_{t_*}^t\int_{B_r(x_*)}\Big\{(|\z_t|+|\z_x|^2)(u-k)_+^2+\z|\z_x|[(u-k)_+]^3+|\r_x|\z^2(u-k)_+
\\\di \qquad\qquad\qquad\qquad+|\r_x|^2\z^2(u-k)_+^2\Big\}dxdt,
\end{array}\label{cut-off}
\end{equation}
where in the last inequality one has used the fact
$$
\begin{array}{ll}
\di |\r^{\a-2}\r_xu_x\z^2(u-k)_+|&\di\leq C|\r_x||\z u_x|\z(u-k)_+\\
&\di =C|\r_x||(\z u)_x-\z_x u|\z(u-k)_+\\
&\di \leq C|\r_x||(\z u)_x|\z(u-k)_++C|\r_x|\z|\z_x|(u-k)_+^2\\
&\di =C|\r_x||[\z (u-k)_+]_x|\z(u-k)_++C|\r_x|\z|\z_x|(u-k)_+^2\\
&\di \leq \f18\r^{\a-1}[\z
(u-k)_+]_x^2+C|\r_x|^2\z^2(u-k)_+^2+C|\z_x|^2(u-k)_+^2.
\end{array}
$$
Thus from \eqref{cut-off}, it holds that
\begin{equation}
\begin{array}{ll}
\di \int_{B_r(x_*)}\z^2(u-k)_+^2(x,t)dx+\int_{t_*}^t\int_{B_r(x_*)}[\z(u-k)_+]_x^2dxdt\\
\di \leq C\int_{B_r(x_*)}\z^2(u-k)_+^2(x,t_*)dx +
C\int_{t_*}^t\int_{B_r(x_*)}\Big\{(|\z_t|+|\z_x|^2)(u-k)_+^2\\
\di \qquad\qquad\qquad\qquad+\z|\z_x|[(u-k)_+]^3+|\r_x|\z^2(u-k)_+
+|\r_x|^2\z^2(u-k)_+^2\Big\}dxdt.
\end{array}\label{cut-off1}
\end{equation}
Now the last three terms in the last integral of \eqref{cut-off1}
can be estimated by
\begin{equation}
\begin{array}{ll}
&\di \int_{t_*}^t\int_{B_r(x_*)}\z|\z_x|[(u-k)_+]^3dxdt \leq
C\int_{t_*}^t\|\z_x(u-k)_+\|_{L^2(B_r(x_*))}\|\z(u-k)_+^2\|_{L^2(B_r(x_*))}dt\\
&\di \leq C\int_{t_*}^t\int_{B_r(x_*)}|\z_x|^2(u-k)_+^2dxdt+C\int_{t_*}^t\int_{B_r(x_*)}\z^2(u-k)_+^4dxdt\\
&\di \leq
C\int_{t_*}^t\int_{B_r(x_*)}|\z_x|^2(u-k)_+^2dxdt+C\int_{t_*}^t\|\z(u-k)_+\|_{L^\i(B_r(x_*))}^2
\|(u-k)_+\|_{L^2(B_r(x_*))}^2dt\\
&\di \leq
C\int_{t_*}^t\int_{B_r(x_*)}|\z_x|^2(u-k)_+^2dxdt+C\int_{t_*}^t\|\z(u-k)_+\|_{L^\i(B_r(x_*))}^2dt\\
&\di \leq C\int_{t_*}^t\int_{B_r(x_*)}|\z_x|^2(u-k)_+^2dxdt
+C\int_{t_*}^t\|\z(u-k)_+\|_{L^2(B_r(x_*))}\|[\z(u-k)_+]_x\|_{L^2(B_r(x_*))}dt\\
&\di \leq C\int_{t_*}^t\int_{B_r(x_*)}|\z_x|^2(u-k)_+^2dxdt
+\f18\int_{t_*}^t\|[\z(u-k)_+]_x\|^2_{L^2(B_r(x_*))}dt\\
&\di\qquad+C\int_{t_*}^t\|\z(u-k)_+\|_{L^2(B_r(x_*))}^2dt,
\end{array}
\label{first}
\end{equation}
\begin{equation}
\begin{array}{ll}
\di \int_{t_*}^t\int_{B_r(x_*)}|\r_x|\z^2(u-k)_+dxdt\\
\di =\iint_{Q_{r,s}^*\cap[u>k]}|\r_x|\z^2(u-k)_+dxdt\\
\di\leq
C\|\z(u-k)_+\|_{L^6(Q_{r,s}^*)}\|\r_x\|_{L^{\f65}(Q_{r,s}^*\cap[u>k])}\\
\di \leq
C\|\z(u-k)_+\|_{V_2(Q_{r,s}^*)}\|\r_x\|_{L^2(Q_{r,s}^*)}|Q_{r,s}^*\cap[u>k]|^{\f13}\\
\di\leq
\f18\|\z(u-k)_+\|^2_{V_2(Q_{r,s}^*)}+C|Q_{r,s}^*\cap[u>k]|^{\f23},
\end{array}\label{second}
\end{equation}
and
\begin{equation}
\begin{array}{ll}
\di \int_{t_*}^t\int_{B_r(x_*)}|\r_x|^2\z^2(u-k)_+^2dxdt\\
\di \leq
C\sup_t\|\r_x\|_{L^2(B_r(x_*))}\int_{t_*}^t\|\z(u-k)_+\|^2_{L^\i(B_r(x_*))}dt\\
\di\leq
C\int_{t_*}^t\|\z(u-k)_+\|_{L^2(B_r(x_*))}\|[\z(u-k)_+]_x\|_{L^2(B_r(x_*))}dt\\
\di \leq
\f18\int_{t_*}^t\|[\z(u-k)_+]_x\|^2_{L^2(B_r(x_*))}dt+C\int_{t_*}^t\|\z(u-k)_+\|_{L^2(B_r(x_*))}^2dt.
\end{array}
\label{third}
\end{equation}
Note that in \eqref{second}, the space $V_2(Q_{r,s}^*)$ is defined
by $V_2(Q_{r,s}^*)=\big\{f\in
L^2(Q_{r,s}^*)~\big|~\|f\|_{V_2(Q_{r,s}^*)}<+\i\big\}, $ with the
norm
$$
\|f\|_{V_2(Q_{r,s}^*)}= {\rm ess}\sup_{\!\!\!\!t\in [t_*,t_*+s]}
\|f(\cdot,t)\|_{L^2(B_r(x_*))}+\|f_x\|_{L^2(Q_{r,s}^*)}.
$$
Substituting \eqref{first}-\eqref{third} into \eqref{cut-off}
implies
\begin{equation}
\begin{array}{ll}
\di \int_{B_r(x_*)}\z^2(u-k)_+^2(x,t)dx+\int_{t_*}^t\int_{B_r(x_*)}[\z(u-k)_+]_x^2dxdt\\
\di \leq\int_{B_r(x_*)}\z^2(u-k)_+^2(x,t_*)dx+C\int_{t_*}^t\int_{B_r(x_*)}\z^2(u-k)_+^2dxdt \\
\di +\f18\|\z(u-k)_+\|^2_{V_2(Q_{r,s}^*)}+
C\int_{t_*}^t\int_{B_r(x_*)}(|\z_t|+|\z_x|^2)(u-k)_+^2dxdt+C|Q_{r,s}^*\cap[u>k]|^{\f23}.
\end{array}\label{E1}
\end{equation}
Applying Gronwall's inequality to \eqref{E1} gives
\begin{equation}
\begin{array}{ll}
\di \sup_{t_*\leq t\leq t_*+s}\int_{B_r(x_*)}\z^2(u-k)_+^2(x,t)dx+\int_{t_*}^{t_*+s}\int_{B_r(x_*)}[\z(u-k)_+]_x^2dxdt\\
\di \leq\int_{B_r(x_*)}\z^2(u-k)_+^2(x,t_*)dx+
C\int_{t_*}^{t_*+s}\int_{B_r(x_*)}(|\z_t|+|\z_x|^2)(u-k)_+^2dxdt+C|Q_{r,s}^*\cap[u>k]|^{\f23}.
\end{array}\label{E2}
\end{equation}
With the  estimate \eqref{E2} at hand, one can show that $u$ is
bounded above locally by the classical De Giorgi-Moser iteration
method and choosing suitably $k$. Similarly, one can obtain the
estimates to $(u-k)_-$ as in \eqref{E2}. Thus we can get the lower
bound for $u$ locally. Furthermore, one can get the local H${\rm
\ddot{o}}$der estimates of $u$ by the classical parabolic theory,
that is, there exists a positive constant $\a_0\in(0,1)$, such that
$$
u\in C^{\a_0,\f{\a_0}{2}}_{\rm loc}(Q^*_{r,s}).
$$

In the following, we will further show that the weak solution
$u(x,t)$ is in fact a strong solution locally as stated in Theorem
\ref{regularity}. Rewrite the momentum equation as
\begin{equation}\label{5.5}
\r u_t+\r uu_x+\g\r^{\g-1}\r_x=\r^\a u_{xx}+\a\r^{\a-1}\r_x u_x.
\end{equation}
Multiplying \eqref{5.5} by $\z^2 u_t$ and integrating the result
over $B_r(x_*)\times (t_*,t]$ with $t\in(t_*,t_*+s)$, one can get
\begin{equation}
\begin{array}{ll}
\di \int_{t_*}^t\int_{B_r(x_*)}\r \z^2u_t^2dxdt\\
\di =\int_{t_*}^t\int_{B_r(x_*)}\Big(-\r u \z^2u_tu_x-\g\r^{\g-1}\r_x\z^2u_t+\r^\a\z^2u_tu_{xx}+\a\r^{\a-1}\r_xu_x\z^2u_t\Big)dxdt\\
\di\leq \f12\int_{t_*}^t\int_{B_r(x_*)}\r
\z^2u_t^2dxdt+C\int_{t_*}^t\int_{B_r(x_*)}\Big(\r
u^2\z^2u_x^2+\r^{2\g-3}\r_x^2\z^2\\
\di\hspace{8cm}+\r^{2\a-1}\z^2u_{xx}^2
+\r^{2\a-3}\z^2\r_x^2u_x^2\Big)dxdt.\\
\end{array}\label{u1}
\end{equation}
By the uniform upper bound and lower bound of the density $\r(x,t)$
in the domain $\Omega_{\s}$ and the local boundedness of the
velocity $u(x,t)$, it follows from \eqref{u1} that
\begin{equation}
\begin{array}{ll}
\di \int_{t_*}^t\int_{B_r(x_*)} \z^2u_t^2dxdt&\di \leq
C\int_{t_*}^t\int_{B_r(x_*)}\Big( \z^2u_x^2+\z^2\r_x^2+\z^2u_{xx}^2
+\z^2\r_x^2u_x^2\Big)dxdt\\
&\di \leq C\int_{t_*}^t\int_{B_r(x_*)}\Big[ \z^2(u-\bar
u)_x^2+\z^2(\r-\bar\r)_x^2+\z^2(\bar u_x^2+\bar\r_x^2)\Big]dxdt\\
&\di \quad+ C\int_{t_*}^t\int_{B_r(x_*)}\Big(\z^2 u_{xx}^2
+\z^2\r_x^2u_x^2\Big)dxdt.
\end{array}\label{ut}
\end{equation}
Note that
\begin{equation}
\begin{array}{ll}
\di \int_{t_*}^t\int_{B_r(x_*)}\z^2\r_x^2u_x^2dxdt&\di \leq
\int_{t_*}^t\|\r_x\|_{L^2(B_r(x_*))}^2\|\z
u_x\|_{L^\i(B_r(x_*))}^2dt\\
&\di \leq
\sup_{t\in(t_*,t_*+s]}\|\r_x\|_{L^2(B_r(x_*))}^2\int_{t_*}^t\|\z
u_x\|_{L^2(B_r(x_*))}\|(\z u_x)_x\|_{L^2(B_r(x_*))}dt\\
&\di \leq \b\int_{t_*}^t\int_{B_r(x_*)}\z^2u_{xx}^2dxdt+C_\b\int_{t_*}^t\int_{B_r(x_*)}(\z^2u_x^2+\z_x^2u_x^2)dxdt\\
&\di \leq \b\int_{t_*}^t\int_{B_r(x_*)}\z^2u_{xx}^2dxdt+C_\b,
\end{array}\label{F1}
\end{equation}
where $\b$ is a small positive constant to be determined and $C_\b$
is the positive constant depending on $\b$.

Thus it follwos from  \eqref{UE},  \eqref{F1} with $\b=1$, and
\eqref{ut} that
\begin{equation}
\int_{t_*}^t\int_{B_r(x_*)} \z^2u_t^2dxdt\leq
C+C\int_{t_*}^t\int_{B_r(x_*)}\z^2u_{xx}^2dxdt.\label{ut1}
\end{equation}
Multiplying  \eqref{vel} by $\z^2 u_{xx}$ and integrating over
$B_r(x_*)\times (t_*,t]$ with $t\in(t_*,t_*+s)$, one can get
\begin{equation}
\begin{array}{ll}
\di \int_{t_*}^t\int_{B_r(x_*)}
\r^{\a-1}\z^2u_{xx}^2dxdt\\
\di=\int_{t_*}^t\int_{B_r(x_*)}\Big(\z^2u_tu_{xx}+\z^2uu_xu_{xx}+\g\r^{\g-2}\r_x\z^2u_{xx}-\a\r^{\a-2}\r_xu_x\z^2u_{xx}\Big)dxdt.
\end{array}
\label{uxx1}
\end{equation} Note that
\begin{equation}
\begin{array}{ll}
\di \z^2u_tu_{xx}&\di =(\z^2u_tu_{x})_x-\z^2u_xu_{xt}-2\z\z_xu_tu_x\\
&\di
=(\z^2u_tu_{x})_x-(\f{\z^2u_x^2}{2})_t+\z\z_tu_x^2-2\z\z_xu_tu_x,
\end{array}
\end{equation}
thus it holds that
\begin{equation}
\int_{t_*}^t\int_{B_r(x_*)}\z^2u_tu_{xx}dxdt=-\int_{B_r(x_*)}\f{\z^2u_x^2}{2}(x,t)dx+\int_{t_*}^t\int_{B_r(x_*)}(\z\z_tu_x^2-2\z\z_xu_tu_x)dxdt.
\label{uxx2}
\end{equation}
Substituting \eqref{uxx2} into \eqref{uxx1} gives
\begin{equation}
\begin{array}{ll}
&\di \f12\int_{B_r(x_*)}\z^2u_x^2(x,t)dx+\int_{t_*}^t\int_{B_r(x_*)}
\z^2u_{xx}^2dxdt\\
&\di
=\int_{t_*}^t\int_{B_r(x_*)}\Big(\z\z_tu_x^2-2\z\z_xu_tu_x+\z^2uu_xu_{xx}+2\r_x\z^2u_{xx}-\f{\r_x}{\r}u_x\z^2u_{xx}\Big)dxdt\\
&\di \leq
\b\int_{t_*}^t\int_{B_r(x_*)}(\z^2u_t^2+\z^2u_{xx}^2)dxdt+C_\b\int_{t_*}^t\int_{B_r(x_*)}\Big[(\z_t^2+\z_x^2)u_x^2+\z^2(\r_x^2+u_x^2)+\z^2\r_x^2u_x^2\Big]dxdt.
\end{array}
\label{uxx3}\end{equation} Combining the estimates \eqref{F1},
\eqref{ut1} and \eqref{uxx3} and choosing both $\b$ suitably small,
we can get
\begin{equation}
\begin{array}{ll}
&\di \f12\int_{B_r(x_*)}\z^2u_x^2(x,t)dx+\int_{t_*}^t\int_{B_r(x_*)}
\z^2(u_t^2+u_{xx}^2)dxdt\leq C. \end{array}
\end{equation}
This completes the proof of Theorem \ref{regularity}.

\section*{Acknowledgments}
Parts of this work were done when the second author Yi Wang was a
postdoctoral fellow in the Institute of Mathematical Science,
Chinese University of Hong Kong during the academic year 2010-2011;
he would like to thank the Institute's support and hospitality.


\begin{thebibliography}{000}


\bibitem{BD1}
D. Bresch and B. Desjardins, Existence of global weak solutions for
a 2D viscous shallow water equations and convergence to the
quasi-geostrophic model, {\it Comm. Math. Phys.}, {\bf 238}(1-2)
(2003), 211-223.

\bibitem{BDL}
D. Bresch, B. Desjardins, Chi-Kun Lin, On some compressible fluid
models: Korteweg, lubrication, and shallow water systems, {\it Comm.
Partial Differential Equations}, {\bf 28}(3-4)(2003), 843-868.

\bibitem{BD2} D. Bresch and B. Desjardins, Some diffusive capillary
models of Korteweg type, {\it C. R. Math. Acad. Sci. Paris, Section
M$\acute{e}$canique, } {\bf 332}(11)(2004), 881-886.

\bibitem{BD3}
D. Bresch and B. Desjardins, On the construction of approximate
solutions for the 2D viscous shallow water model and for
compressible Navier-Stokes models, {\it J. Math. Pures Appl.}, {\bf
86} (2006), 362-368.




\bibitem{BDG}
D. Bresch, B. Desjardins, D. Gerard-Varet, On compressible
Navier-Stokes equations with density dependent viscosities in
bounded domains, {\it J. Math. Pures Appl.}, {\bf 87}(2) (2007),
227-235.

\bibitem{Chen}
G. Q. Chen, Vacuum states and global stability of rarefaction waves
for compressible flow, {\it  Methods Appl. Anal.}, {\bf 7} (2000),
337-361.

\bibitem{FZ} D. Y. Fang, T. Zhang, Compressible Navier-Stokes
equations with vacuum state in the case of general pressure law,
{\it Math. Methods Appl. Sci.}, {\bf 29} (2006), 1081-1106.

\bibitem{FNP00} E. Feireisl, A. Novotn\'y and H. Petzeltov\'a,
On the existence of globally defined weak solutions to the
Navier-Stokes equations of isentropic compressible fluids, {\it J.
Math. Fluid Mech.}, {\bf 3} (2001), 358-392.

\bibitem{GP} J. F. Gerbeau, B. Perthame, Derivation of viscous
Saint-Venant system for laminar shallow water; numerical validation,
{\it Discrete Contin. Dyn. Syst. (Ser. B)}, {\bf 1} (2001), 89-102.

\bibitem{GJX} Z. H. Guo, Q. S. Jiu, Z. P. Xin, Spherically
symmetric isentropic compressible flows with density-dependent
viscosity coefficients,  {\it SIAM J. Math. Anal.}, {\bf 39}(5)
(2008), 1402-1427.


\bibitem{Hoff1} D. Hoff, Global existence for 1D, compressible,
isentropic Navier-Stokes equations with large initial data, {\it
Trans. Amer. Math. Soc.}, {\bf 303}(1)(1987), 169-181.

\bibitem{Hoff} D. Hoff, Strong convergence to global solutions for
multidimensional flows of compressible, viscous fluids with
polytropic equations of state and discontinuous initial data, {\it
Arch. Rat. Mech. Anal.}, {\bf 132}(1995), 1-14.

\bibitem{HD} D. Hoff, D. Serre, The failure of continuous
dependence on initial data for the Navier-Stokes equations of
compressible flow, {\it SIAM J. Appl. Math.}, {\bf 51}(1991),
887-898.

\bibitem{Hoff-Smoller} D. Hoff, J. Smoller,
Non-formation of vacuum states for compressible Navier-Stokes
equations, {\it  Comm. Math. Phys.}, {\bf 216}(2) (2001), 255-276.

\bibitem{HLW}F. Huang,  M. Li, Y. Wang, Zero dissipation limit to rarefaction waves with vacuum for 1-D
compressible isentropic Navier-Stokes equations, Preprint.

\bibitem{J}
S. Jiang, Global smooth solutions of the equations of a viscous,
heat-conducting one-dimensional gas with density-dependent
viscosity, {\it Math. Nachr.}, {\bf 190}(1998), 169-183.


\bibitem{JXZ} S. Jiang, Z. P. Xin and P. Zhang, Global weak
solutions to 1D compressible isentropy Navier-Stokes with
density-dependent viscosity,  {\it Methods and Applications of
Analysis}, {\bf 12} (3)(2005), 239-252.

\bibitem{JWX} Q. S. Jiu, Y. Wang, Z. P. Xin, Stability of rarefaction waves
 to the 1D compressible Navier-Stokes equations with density-dependent
 viscosity,
 {\it Comm. Part. Diff. Equ.}, {\bf 36} (2011), 602-634.

\bibitem{JX} Q. S. Jiu, Z. P. Xin, The Cauchy problem for
 1D compressible flows with density-dependent viscosity
 coefficients,
 {\it Kinet. Relat. Models},  {\bf 1} (2),  (2008),   313-330.

 \bibitem{Ka}
J. I. Kanel, A model system of equations for the one-dimensional
motion of a gas, {\it Diff. Uravn.}, {\bf 4}, (1968), 721-734 (in
Russian).

\bibitem{KS}
A. V. Kazhikhov, V. V. Shelukhin, Unique global solution with
respect to time of initial-boundary value problems for
one-dimensional equations of a viscous gas, {\it J. Appl. Math.
Mech. }{\bf 41}(1977), 273-282; translated from {\it Prikl. Mat.
Meh. }{\bf 41 }(1977), 282-291.

\bibitem{LLX}
H. L. Li, J. Li, Z. P. Xin, Vanishing of vacuum states and blow-up
phenomena of the compressible Navier-Stokes equations, {\it Comm.
Math. Phys.}, {\bf 281}(2),  (2008),  401--444.

\bibitem{Li-J}
W. M. Li, Q. S. Jiu, Global weak solutions to one-dimensional
compressible Navier-Stokes equations with density-dependent
viscosity coefficients, {\it J. Partial Differ. Equ.},  {\bf 23} (3)
(2010), 290-304.

\bibitem{L98}
P. L. Lions, {\it Mathematical Topics in Fluid Dynamics 2,
Compressible Models}, Oxford Science Publication, Oxford, 1998.

\bibitem{lius}
T.-P. Liu, J. Smoller,  On the vacuum state for the isentropic
 gas
dynamics equations, {\it Adv. in Appl. Math.},
 {\bf1}(4) (1980), 345-359.

\bibitem{Liu-Xin-1} T. P. Liu, Z. P. Xin, Nonlinear stability of rarefaction
waves for compressible Navier-Stokes equations, {\it Comm. Math.
Phys.}, {\bf118} (1988), 451--465.


\bibitem{LXY}
T. P. Liu, Z. P. Xin, T. Yang, Vacuum states of compressible flow,
{\it Discrete Continuous Dynam. systems}, {\bf 4}(1)(1998), 1-32.

\bibitem{LuoXY}
T. Luo, Z. P.  Xin, T. Yang, Interface behavior of compressible
Navier-Stokes equations with vacuum, {\it  SIAM J. Math. Anal.},
{\bf 31} (2000), 1175-1191.


\bibitem{Matsumura-Nishihara-1} A. Matsumura, K. Nishihara, Asymptotics toward
the rarefaction wave of the solutions of a one-dimensional model
system for compressible viscous gas, {\it Japan J. Appl. Math.},
{\bf3} (1986), 1--13.

\bibitem{Matsumura-Nishihara-2} A. Matsumura, K. Nishihara, Global stability of the
rarefaction wave of a one-dimensional model system for compressible
viscous gas,  {\it Comm. Math. Phys.},  {\bf 144}(2),  (1992),
325--335.

\bibitem{MV}
A. Mellet and A. Vasseur, On the barotropic compressible
Navier-Stokes equation, {\it Comm. Partial Differential Equations},
{\bf 32}(3) (2007), 431-452.

\bibitem{MV2}
A. Mellet and A. Vasseur, Existence and uniqueness of global strong
solutions for one-dimensional compressible Navier-Stokes equations,
{\it SIAM J. Math. Anal.}, {\bf 39}(4) (2008), 1344-1365.

\bibitem{OMM} M. Okada, $\breve{S}$.
Matu$\breve{s}$u-Ne$\breve{c}$asov$\acute{a}$, T. Makino, Free
boundary problem for the equation of one-dimensional motion of
compressible gas with density-dependent viscosity, {\it Ann. Univ.
Ferrara Sez. VII (N.S.)}, {\bf 48}(2002), 1-20.

\bibitem{R} O. Rozanova, Blow up of smooth solutions to the compressible Navier-Stokes
equations with the data highly decreasing at infinity,
 {\it J. Differ. Eqs.}, \textbf{245} (2008), 1762-1774.

\bibitem{SZ} I. Straskraba, A. Zlotnik, Global properties of
solutions to 1D viscous compressible barotropic fluid equations with
density dependent viscosity, {\it Z. Angew. Math. Phys.}, {\bf 54
}(4) (2003), 593-607.

\bibitem{VYZ}
S. W. Vong, T. Yang, C. J. Zhu, Compressible Navier-Stokes equations
with degenerate viscosity coefficient and vacuum II, {\it J.
Differential Equations}, {\bf 192}(2)(2003), 475-501.

\bibitem{Xin}

Z. P. Xin, Blow-up of smooth solution to the compressible
Navier-Stokes equations with compact density, {\it Comm. Pure  Appl.
Math.},  {\bf 51}(1998), 229-240.

\bibitem{YYZ}
T. Yang, Z. A. Yao, C. J. Zhu, Compressible Navier-Stokes equations
with density-dependent viscosity and vacuum, {\it Comm. Partial
Differential Equations}, {\bf 26} (5-6)(2001), 965-981.



\bibitem{YZ2}
T. Yang, C. J. Zhu, Compressible Navier-Stokes equations with
degenerate viscosity coefficient and vacuum, {\it Comm. Math.
Phys.}, {\bf 230} (2)(2002), 329-363.


 \end{thebibliography}
\end{document}